\documentclass[10pt,dvips,twoside]{article}

\usepackage{amssymb,amsmath,makeidx}
\usepackage[T1]{fontenc}
\usepackage{latexsym}
\usepackage{exscale}
\usepackage{pb-diagram}
\usepackage[dvips]{graphicx}
\usepackage{xr}

\newcommand{\al}{\alpha}
\newcommand{\alvec}{{\vec{\al}}} 
\newcommand{\alvecbe}{{\vec{\al}^\frown\be}} 
\newcommand{\alvecga}{{\vec{\al}^\frown\ga}} 
\newcommand{\alvecpr}{{\vec{\al}'}}
\newcommand{\alvecpl}{{\vec{\al}^+}}

\newcommand{\alpr}{{\alpha^\prime}}
\newcommand{\alnod}{{\al_0}}

\newcommand{\ale}{{\al_1}}

\newcommand{\ali}{{\al_i}}
\newcommand{\alinod}{{\al_{i_0}}}

\newcommand{\alie}{{\al_{i+1}}}

\newcommand{\alimin}{{\al_{i-1}}}

\newcommand{\alj}{{\al_j}}
\newcommand{\alje}{{\al_{j+1}}}
\newcommand{\aljnod}{{\al_{j_0}}}

\newcommand{\alk}{{\al_k}}

\newcommand{\aln}{{\al_n}}
\newcommand{\alne}{{\al_{n+1}}}

\newcommand{\alnmin}{{\al_{n-1}}}

\newcommand{\alvecrestrnmin}{{\alvec{\scriptscriptstyle{\restriction {n-1}}}}}

\newcommand{\albar}{\bar{\al}}
\newcommand{\alnbar}{\overline{\aln}}

\newcommand{\alhat}{{\widehat{\al}}}
\newcommand{\be}{\beta}
\newcommand{\beti}{\tilde{\be}}
\newcommand{\bevec}{{\vec{\be}}}

\newcommand{\bebar}{\bar{\be}}
\newcommand{\bebarpr}{\bar{\be}^\prime}
\newcommand{\bepr}{{\beta^\prime}}

\newcommand{\behat}{{\widehat{\be}}}
\newcommand{\ga}{\gamma}
\newcommand{\gahat}{{\widehat{\ga}}}
\newcommand{\gapr}{{\ga^\prime}}
\newcommand{\gati}{\tilde{\gamma}}

\newcommand{\gavec}{{\vec{\ga}}}

\newcommand{\de}{\delta}

\newcommand{\deti}{\tilde{\delta}}

\newcommand{\devec}{{\vec{\de}}}

\newcommand{\devecpr}{{\vec{\de}^\prime}}
\newcommand{\De}{\Delta}

\newcommand{\eps}{\varepsilon}

\newcommand{\la}{\lambda}
\newcommand{\lapr}{{\la^\prime}}
\newcommand{\ka}{\kappa}
\newcommand{\om}{\omega}

\newcommand{\Om}{\Omega}

\newcommand{\sivec}{{\vec{\si}}}

\newcommand{\sipr}{{\si^\prime}}

\newcommand{\tauti}{\tilde{\tau}}
\newcommand{\tauvec}{{\vec{\tau}}}
\newcommand{\tauvecsi}{{\vec{\tau}^\frown\si}}

\newcommand{\taupr}{{\tau^\prime}}

\newcommand{\taunpr}{{\tau^\prime_n}}

\newcommand{\taui}{{\tau_i}}
\newcommand{\tauie}{{\tau_{i+1}}}
\newcommand{\taujstar}{{\tau_j^\star}}

\newcommand{\taunstar}{{\tau_n^\star}}

\newcommand{\si}{\sigma}
\newcommand{\tht}{\vartheta}

\newcommand{\ze}{\zeta}

\newcommand{\etavec}{{\vec{\eta}}}

\renewcommand{\phi}{\varphi}

\newcommand{\N}{{\mathbb N}}
\newcommand{\Hz}{{\mathbb P}}
\newcommand{\Lz}{{\mathbb L}}
\newcommand{\Ez}{{\mathbb E}}
\newcommand{\Ezone}{{\mathbb E}_1}
\newcommand{\On}{{\mathrm{Ord}}}

\newcommand{\CNF}{{\mathrm{\scriptscriptstyle{CNF}}}}
\newcommand{\ANF}{{\mathrm{\scriptscriptstyle{ANF}}}}

\newcommand{\NF}{{\mathrm{\scriptscriptstyle{NF}}}}

\newcommand{\Lim}{\mathrm{Lim}}
\newcommand{\Image}{\mathrm{Im}}

\newcommand{\finsub}{\subseteq_\mathrm{fin}}
\newcommand{\logend}{{\mathrm{logend}}}
\newcommand{\sumend}{{\mathrm{end}}}

\newcommand{\leo}{\le_1}

\newcommand{\lo}{<_1}

\newcommand{\klex}{<_\mathrm{\scriptscriptstyle{lex}}}
\newcommand{\kglex}{\le_\mathrm{\scriptscriptstyle{lex}}}

\newcommand{\thtk}{\tht_k}

\newcommand{\thtt}{\tht^\tau}
\newcommand{\thtti}{\tht^\taui}

\newcommand{\Targ}[1]{{\operatorname{T}^{#1}}}

\newcommand{\Tt}{{\operatorname{T}^\tau}}
\newcommand{\Ttn}{{\operatorname{T}^{\tau_n}}}
\newcommand{\Tti}{{\operatorname{T}^{\taui}}}
\newcommand{\Ttvec}{{\operatorname{T}^\tauvec}}

\newcommand{\ltvec}{{\operatorname{l}^\tauvec}}

\newcommand{\ltvecbe}{{\operatorname{l}^{\tauvec^\frown\be}}}
\newcommand{\lalvec}{{\operatorname{l}^\alvec}}

\newcommand{\htarg}[1]{\operatorname{ht}_{#1}}
\newcommand{\hte}{{\operatorname{ht}_1}}

\newcommand{\lh}{{\operatorname{lh}}}
\newcommand{\gs}{{\operatorname{gs}}}

\newcommand{\kval}{\kappa^\alvec}

\newcommand{\kvalnminaln}{{\kappa^\alvecrestrnmin_\aln}}

\newcommand{\kvalbe}{\kappa^\alvec_\be}

\newcommand{\kvalga}{\kappa^\alvec_\ga}

\newcommand{\kvalde}{\kappa^\alvec_\de}

\newcommand{\kvga}{\kappa^\gavec}

\newcommand{\laaln}{{\la_{\aln}}}
\newcommand{\laalnbe}{{\la^\aln_\be}}

\newcommand{\iotal}{\iota_{\tau,\al}}

\newcommand{\zetal}{{\ze^\tau_\al}}
\newcommand{\zetbe}{{\ze^\tau_\be}}

\newcommand{\latal}{{\la^\tau_\al}}

\newcommand{\latbe}{{\la^\tau_\be}}

\newcommand{\pist}{\pi_{\si,\tau}}

\newcommand{\mc}{{\operatorname{mc}}}
\newcommand{\mf}{{\operatorname{mf}}}
\newcommand{\lf}{{\operatorname{lf}}}

\newcommand{\kplnod}{{\operatorname{KP}\!\ell_0}}

\newcommand{\Rtwo}{{{\cal R}_2}}
\newcommand{\Ctwo}{{{\cal C}_2}}

\newcommand{\Ronepl}{{{\cal R}_1^+}}
\newcommand{\Rtwopl}{{{\cal R}_2^+}}

\newcommand{\bardot}{\bar{\cdot}}

\newcommand{\qed}{\mbox{ }\hfill $\Box$\vspace{2ex}}
\newcommand{\imp}{\Rightarrow}
\newcommand{\aeq}{\Leftrightarrow}
\newcommand{\andsp}{\:\&\:}

\newcommand{\sub}{\subseteq}

\newcommand{\set}[2]{\{ #1 \:|\: #2\}}

\newcommand{\singleton}[1]{\{ #1 \}}

\newlength{\hilflh}
\newcommand{\hilfminus}[1]{
  \settowidth{\hilflh}{$#1-$}\mbox{$#1-\hspace{-0.5\hilflh}
  \makebox[0pt]{\raisebox{0.24\hilflh}{$#1\cdot$}}\hspace{0.5\hilflh}$}}
\newcommand{\minusp}{\mathbin{\mathchoice {\hilfminus{\displaystyle}}
  {\hilfminus{\textstyle}}{\hilfminus{\scriptstyle}}
  {\hilfminus{\scriptscriptstyle}}}}

\newtheorem{theo}{Theorem}[section]
\newtheorem{cor}[theo]{Corollary}
\newtheorem{lem}[theo]{Lemma}
\newtheorem{defi}[theo]{Definition}
\newtheorem{prop}[theo]{Proposition}

\newcommand{\oneinf}{1^\infty}

\newcommand{\tauinf}{\tau^\infty}
\newcommand{\alinf}{\al^\infty}

\newcommand{\chial}{\chi^\al}

\newcommand{\chialncheck}{\check{\chi}^\aln}

\newcommand{\mutal}{\mu^\tau_\al}
\newcommand{\mual}{\mu_\al}
\newcommand{\mualj}{\mu_\alj}
\newcommand{\mualn}{\mu_\aln}

\newcommand{\mutali}{\mu^\tau_\ali}

\newcommand{\mutaln}{\mu^\tau_\aln}

\newcommand{\mube}{\mu_\be}
\newcommand{\muga}{\mu_\ga}

\newcommand{\rhoalmutal}{\varrho^\al_{\mutal}}
\newcommand{\rhoalbe}{{\varrho^\al_\be}}

\newcommand{\rhoalnga}{{\varrho^\aln_\ga}}

\newcommand{\MNF}{{\mathrm{\scriptscriptstyle{MNF}}}}
\newcommand{\mNF}{{\mathrm{\scriptscriptstyle{NF}}}}
\newcommand{\Mz}{{\mathbb M}}

\newcommand{\trs}{{\mathrm{ts}}}
\newcommand{\trst}{{\mathrm{ts}^\tau}}
\newcommand{\trsal}{{\mathrm{ts}^\al}}

\newcommand{\trsaln}{{\mathrm{ts}^\aln}}

\newcommand{\trsbe}{{\mathrm{ts}^\be}}

\newcommand{\cs}{{\mathrm{cs}}}
\newcommand{\cspr}{{\mathrm{cs}^\prime}}

\newcommand{\tc}{{\mathrm{tc}}}
\newcommand{\TC}{{\mathrm{TC}}}

\newcommand{\RS}{\mathrm{RS}}
\newcommand{\RSt}{\mathrm{RS}^\tau}
\newcommand{\letwo}{\le_2}
\newcommand{\ktwo}{<_2}

\newcommand{\lSeq}{\mathrm{lSeq}}

\newcommand{\dom}{\mathrm{dom}}

\newcommand{\domkval}{{\mathrm{dom}({\kval})}}

\newcommand{\domnuval}{{\mathrm{dom}({\nuval})}}

\newcommand{\dpf}{\mathrm{dp}}

\newcommand{\dpval}{{\mathrm{dp}_\alvec}}

\newcommand{\dpvga}{{\mathrm{dp}_\gavec}}

\newcommand{\nuval}{\nu^\alvec}

\newcommand{\nuvalbe}{\nu^\alvec_\be}
\newcommand{\nuvalga}{\nu^\alvec_\ga}

\newcommand{\nuvga}{\nu^\gavec}

\newcommand{\alcp}[1]{{\al_{#1}}}
\newcommand{\becp}[1]{{\be_{#1}}}

\newcommand{\taucp}[1]{{\tau_{#1}}}
\newcommand{\taucppr}[1]{{\tau^\prime_{#1}}}

\newcommand{\tauticp}[1]{{\tauti_{#1}}}

\newcommand{\ov}{\mathrm{o}}

\newcommand{\ordcp}[1]{{\mathrm{o}_{#1}}}

\newcommand{\letc}{\le_\mathrm{TC}}

\newcommand{\cml}{\operatorname{cml}}
\newcommand{\gbo}{\operatorname{gbo}}
\newcommand{\predec}{\operatorname{pred}}
\newcommand{\predecs}{\operatorname{Pred}}

\newcommand{\TS}{\operatorname{TS}}
\newcommand{\TSe}{\operatorname{TS}^1}
\newcommand{\TSt}{\operatorname{TS}^\tau}

\newcommand{\TSal}{\operatorname{TS}^\al}

\newcommand{\mts}{\operatorname{mts}}
\newcommand{\mtsal}{\operatorname{mts}^\al}
\newcommand{\mtsale}{\operatorname{mts}^\ale}
\newcommand{\mtsbe}{\operatorname{mts}^\be}

\newcommand{\hop}{\operatorname{h}}
\newcommand{\homega}{\operatorname{h}_\om}
\newcommand{\hbe}{\operatorname{h}_\be}
\newcommand{\hga}{\operatorname{h}_\ga}
\newcommand{\hde}{\operatorname{h}_\de}

\newcommand{\maxmucov}{\operatorname{max-cov}}
\newcommand{\maxmucovtau}{\operatorname{max-cov}^\tau}
\newcommand{\minmucoval}{\operatorname{min-cov}^\al}

\newcommand{\sk}{\operatorname{sk}}
\newcommand{\skbe}{\operatorname{sk}_\be}
\newcommand{\skga}{\operatorname{sk}_\ga}

\newcommand{\me}{\operatorname{me}}

\def\vec#1{\mathchoice{\mbox{\boldmath$\displaystyle#1$}}
{\mbox{\boldmath$\textstyle#1$}}
{\mbox{\boldmath$\scriptstyle#1$}}
{\mbox{\boldmath$\scriptscriptstyle#1$}}}

\begin{document}

\title{Tracking chains revisited}

\author{Gunnar Wilken
\footnote{This article is a pre-print of a chapter in {\it Sets and Computations}, Lecture Notes Series Vol.\ 33, 
Institute for Mathematical Sciences, National University of Singapore, \copyright World Scientific Publishing Company (2017), see \cite{W17}.
The author would like to acknowledge the Institute for Mathematical Sciences of the National University of Singapore
for its partial support of this work during the ``Interactions'' week of the workshop {\it Sets and Computations} in April 2015.
}\\
Structural Cellular Biology Unit\\
Okinawa Institute of Science and Technology\\
1919-1 Tancha, Onna-son, 904-0495 Okinawa, Japan\\
{\tt wilken@oist.jp}
}

\maketitle

\begin{abstract}
The structure $\Ctwo:=(1^\infty,\le,\leo,\letwo)$, introduced and first analyzed in \cite{CWc}, is shown to be elementary recursive.
Here, $1^\infty$ denotes the proof-theoretic ordinal of the fragment $\Pi^1_1$-$\mathrm{CA}_0$ of second order number theory, 
or equivalently the set theory $\kplnod$, which axiomatizes limits of models of Kripke-Platek set theory with infinity.
The partial orderings $\leo$ and $\letwo$ denote the relations of $\Sigma_1$- and $\Sigma_2$-elementary substructure, respectively.
In a subsequent article \cite{W} we will show that the structure $\Ctwo$ comprises the core of the structure $\Rtwo$ of pure 
elementary patterns of resemblance of order $2$.
In \cite{CWc} the stage has been set by showing that the least ordinal containing a cover of each pure pattern of order $2$
is $1^\infty$. However, it is not obvious from \cite{CWc} that $\Ctwo$ is an elementary recursive structure. This is shown here
through a considerable disentanglement in the description of connectivity components of $\leo$ and $\letwo$.
The key to and starting point of our analysis is the apparatus of ordinal arithmetic developed in \cite{W07a} and in Section $5$ 
of \cite{CWa}, which was enhanced in \cite{CWc}, specifically for the analysis of $\Ctwo$.
\end{abstract}

\section{Introduction}
Let $\Rtwo=\left(\On;\le,\leo,\letwo\right)$ be the structure of ordinals with standard linear ordering $\le$ and 
partial orderings $\le_1$ and $\le_2$, simultaneously defined by induction on $\be$ in  
\[\al\le_i\be:\aeq \left(\al;\le,\le_1,\le_2\right) \preceq_{\Sigma_i} \left(\be;\le,\le_1,\le_2\right)\]
where $\preceq_{\Sigma_i}$ is the usual notion of $\Sigma_i$-elementary substructure (without bounded quantification), see \cite{C99,C01}
for fundamentals and groundwork on elementary patterns of resemblance.
Pure patterns of order $2$ are the finite isomorphism types of $\Rtwo$. The \emph{core}
of $\Rtwo$ consists of the union of \emph{isominimal realizations} of these patterns within $\Rtwo$, where a finite 
substructure of $\Rtwo$ is called isominimal, if it is pointwise minimal (with respect to increasing enumerations) 
among all substructures of $\Rtwo$ isomorphic to it, and where an isominimal substructure of $\Rtwo$ realizes a pattern $P$, if
it is isomorphic to $P$. It is a basic observation, cf.\ \cite{C01}, that the class of pure patterns of order $2$ is contained in the class 
$\mathcal{RF}_2$ of \emph{respecting forests of order $2$}:
finite structures $P$ over the language $(\le_0,\leo,\letwo)$ where $\le_0$ is a linear ordering and $\leo,\letwo$ are forests such that 
$\letwo\subseteq\leo\subseteq\le_0$ and $\le_{i+1}$ \emph{respects} $\le_i$, i.e.\ $p\le_i q\le_i r\andsp p\le_{i+1}r$ implies $p\le_{i+1}q$
for all $p,q,r\in P$, for $i=0,1$.  

In \cite{CWc} we showed that every pattern has a cover below $1^\infty$, the least such ordinal.
Here, an order isomorphism (embedding) is a cover (covering, respectively) if it maintains the relations $\le_1$ and $\le_2$.
The ordinal of $\kplnod$ is therefore least such that there exist arbitrarily long finite $\le_2$-chains.
Moreover, by determination of enumeration functions of (relativized) connectivity components of $\le_1$ and $\le_2$ we
were able to describe these relations in terms of classical ordinal notations. The central observation in connection with
this is that every ordinal below $1^\infty$ is the greatest element in a $\le_1$-chain in which $\le_1$- and $\le_2$-chains alternate.
We called such chains \emph{tracking chains} as they provide all $\le_2$-predecessors and the greatest $\le_1$-predecessors
insofar as they exist. 

In the present article we will review and slightly extend the ordinal arithmetical toolkit and then 
verify through a disentangling reformulation, that \cite{CWc} in fact yields an elementary recursive characterization of the restriction of $\Rtwo$ to the structure 
$\Ctwo=(1^\infty;\le,\le_1,\le_2)$. It is not obvious from \cite{CWc} that $\Ctwo$ is an elementary recursive structure since several proofs 
there make use of transfinite induction up to $1^\infty$, which allowed for a somewhat shorter argumentation there.  
We will summarize the results in \cite{CWc} to a sufficient and convenient degree.
As a byproduct, \cite{CWc} will become considerably more accessible. We will prove the equivalence of the arithmetical descriptions of 
$\Ctwo$ given in \cite{CWc} and here.
Note that the equivalence of this elementary recursive characterization with the original structure based on elementary substructurehood is proven in Section $7$ of \cite{CWc}, using full transfinite induction up to the 
ordinal of $\kplnod$. In this article we rely on this result and henceforth identify $\Ctwo$ with its arithmetical characterization given in \cite{CWc} and further illuminated in Section \ref{revisitsec} of the present article, where we also show that the finite isomorphism types
of the arithmetical $\Ctwo$ are respecting forests of order $2$, without relying on semantical characterization of the arithmetical $\Ctwo$.

With these preparations out of the way we will be able to provide, in a subsequent article \cite{W}, an algorithm that assigns 
an isominimal realization within $\Ctwo$ to each respecting forest
of order $2$, thereby showing that each such respecting forest is in fact (up to isomorphism)
a pure pattern of order $2$.
The approach is to formulate the corresponding theorem flexibly so that isominimal realizations above certain relativizing
tracking chains are considered.   
There we will also define an elementary recursive function that assigns descriptive patterns $P(\al)$ to
ordinals $\al\in 1^\infty$. A descriptive pattern for an ordinal $\al$ in the above sense is a pattern, the isominimal 
realization of which contains $\al$. 
Descriptive patterns will be given in a way that makes a canonical choice for normal
forms, since in contrast to the situation in $\Ronepl$, cf.\ \cite{W07c,CWa}, there is no unique notion of normal form in $\Rtwo$.
The chosen normal forms will be of least possible cardinality.

The mutual order isomorphisms between hull and pattern notations that will be given in \cite{W} enable classification of a new independence
result for $\kplnod$: We will demonstrate that the result by Carlson in \cite{C16}, according to which the collection of respecting forests of 
order $2$ is well-quasi ordered with respect to coverings, 
cannot be proven in $\kplnod$ or, equivalently, in the restriction $\Pi^1_1\mathrm{-CA}_0$ of second order number theory to 
$\Pi^1_1$-comprehension and set induction. On the other hand,  
we know that transfinite induction up to the ordinal $1^\infty$ of $\kplnod$ suffices to show that every pattern is covered \cite{CWc}.

This article therefore delivers the first part of an in depth treatment of the insights and results presented in a lecture during 
the ``Interactions'' week
of the workshop {\it Sets and Computations}, held at the Institute for Mathematical Sciences of the National University of Singapore in April 2015.

\section{Preliminaries}
The reader is assumed to be familiar with basics of ordinal arithmetic (see e.g.\ \cite{P09})
and the ordinal arithmetical tools developed in \cite{W07a} and Section 5 of \cite{CWa}. See the index at the
end of \cite{W07a} for quick access to its terminology. Section 2 of \cite{W07c} (2.1--2.3) provides a summary of results
from \cite{W07a}. As mentioned before, we will build upon \cite{CWc}, the central concepts and results of which will be reviewed 
here and in the next section. For detailed reference see also the index of \cite{CWc}.

\subsection{Basics}
Here we recall terminology already used in \cite{CWc} (Section 2) for the reader's convenience.
Let $\Hz$ denote the class of additive principal numbers, i.e.\ nonzero ordinals that are closed under ordinal addition,
that is the image of ordinal exponentiation to base $\om$.
Let $\Lz$ denote the class of limits of additive principal numbers, i.e.\ the limit points of $\Hz$, 
and let $\Mz$ denote the class of multiplicative principal numbers, i.e.\ nonzero ordinals closed under ordinal multiplication.
By $\Ez$ we denote the class of epsilon numbers, i.e.\ the fixed-points of $\om$-exponentiation.

We write $\al=_\ANF\ale+\ldots+\aln$ if $\ale,\ldots,\aln\in\Hz$ such that $\ale\ge\ldots\ge\aln$, which is called the 
representation of $\al$ in additive normal form,
and $\al=_\NF\be+\ga$ if the expansion of $\be$ into its additive normal form (ANF) in the sum $\be+\ga$ syntactically results in the 
additive normal form of $\al$.
The Cantor normal form representation of an ordinal $\al$ is given by $\al=_\CNF\om^\ale+\ldots+\om^\aln$ where $\ale\ge\ldots\ge\aln$
with $\al>\ale$ unless $\al\in\Ez$.
For $\al=_\ANF\ale+\ldots+\aln$, we define $\mc(\al):=\ale$\index{$\mc$} and $\sumend(\al):=\aln$.
We set $\sumend(0):=0$.
Given ordinals $\al, \be$ with $\al\le\be$ we write $-\al+\be$\index{$-\al+\be$} for the unique $\ga$ such that $\al+\ga=\be$.
As usual let  $\al\minusp\beta$\index{$-p$@$\al\minusp\be$} be $0$ if $\be\ge\al$, $\ga$ if $\be<\al$ and there exists the minimal $\ga$ s.t.\  $\al=\ga+\be$, and $\al$ otherwise.

For $\al\in\On$ we denote the least multiplicative principal number greater than $\al$ by $\al^\Mz$.
Notice that if $\al\in\Hz$, $\al>1$, say $\al=\om^\alpr$, we have $\al^\Mz=\al^\om=\om^{\alpr\cdot\om}$.
For $\al\in\Hz$ we use the following notations for \emph{multiplicative normal form}:\index{multiplicative normal form}
\begin{enumerate}
\item $\al=_\mNF\eta\cdot\xi$\index{$=_\mNF$} if and only if $\xi=\om^{\xi_0}\in\Mz$ (i.e.\ $\xi_0\in\singleton{0}\cup\Hz$) and
either $\eta=1$ or $\eta=\om^{\eta_1+\ldots+\eta_n}$ such that $\eta_1+\ldots+\eta_n+\xi_0$ is in
additive normal form. When ambiguity is unlikely, we sometimes allow $\xi$ to be of a form $\om^{\xi_1+\ldots+\xi_m}$ such that 
$\eta_1+\ldots+\eta_n+\xi_1+\ldots+\xi_m$ is in additive normal form.
\item $\al=_\MNF\ale\cdot\ldots\cdot\alk$\index{$=_\MNF$} if and only if $\ale,\ldots,\alk$ is the unique decreasing sequence of multiplicative principal numbers, the product of which is equal to $\al$.
\end{enumerate}
For $\al\in\Hz$, $\al=_\MNF\ale\cdot\ldots\cdot\alk$, we write $\mf(\al)$ for $\ale$ and $\lf(\al)$\index{$\lf$} for $\alk$. 
Note that if $\al\in\Hz-\Mz$ then $\lf(\al)\in\Mz^{>1}$ and $\al=_\mNF\albar\cdot\lf(\al)$
where the definition of $\albar$ given in \cite{W07a} for limits of additive principal numbers is
extended to ordinals $\al$ of a form $\al=\om^{\alpr+1}$ by $\albar:=\om^\alpr$, see Section 5 of \cite{CWa}.

Given $\al,\be\in\Hz$ with $\al\le\be$ we write $(1/\al)\cdot\be$\index{$/$@$(1/\al)\cdot\be$} for the uniquely determined ordinal $\ga\le\be$
such that $\al\cdot\ga=\be$. Note that with the representations $\al=\om^\alpr$ and $\be=\om^\bepr$ we have
\[(1/\al)\cdot\be=\om^{-\alpr+\bepr}.\]

For any $\al$ of a form $\om^\alpr$ we write $\log(\al)$\index{$\log$} for $\alpr$, and we set $\log(0):=0$. For an arbitrary ordinal $\be$ we
write $\logend(\be)$ for $\log(\sumend(\be))$.

\subsection{Relativized notation systems $\Tt$}
Settings of relativization are given by ordinals from $\Ezone:=\singleton{1}\cup\Ez$\index{$\ezo$@$\Ezone$} and frequently indicated by Greek letters, preferably $\si$ or $\tau$.
Clearly, in this context $\tau=1$ denotes the trivial setting of relativization. For a setting $\tau$ of relativization we define 
$\tauinf:=\Tt\cap\Omega_1$\index{$\tauinf$}
where $\Tt$ is defined as in \cite{W07a} and reviewed in Section 2.2 of \cite{W07c}.
$\Tt$ is the closure of parameters below $\tau$ under addition and the stepwise, injective, and fixed-point free collapsing functions 
$\thtk$ the domain of which is $\Tt\cap\Om_{k+2}$, where $\thtt:=\tht_0$ is relativized to $\tau$. 
As in \cite{CWc}, most considerations will be confined to the segment $\oneinf$.
Translation between different settings of relativization, see Section 6 of \cite{W07a}, is effective on the term syntax and enjoys convenient invariance properties regarding the operators described below, as was verified in \cite{W07a}, \cite{CWa}, and \cite{CWc}. 
We therefore omit the purely technical details here. 

\subsection{Refined localization}\label{reflocsubsec}
The notion of $\tau$-localization (Definition 2.11 of \cite{W07c}) and its refinement to $\tau$-fine-localization by iteration of
the operator $\bardot$, see Definitions 5.1 and 5.5 of \cite{CWa}, continue to be essential as they locate ordinals in terms of 
closure properties (fixed-point level and limit point thinning). These notions are effectively derived from the term syntax.
We refer to Subsection 2.3 of \cite{W07c} and Section 5 of \cite{CWa} for a complete picture of these concepts. In the present article,
the operator $\bardot$ is mostly used to decompose ordinals that are not multiplicative principal, i.e. if $\al=_\mNF\eta\cdot\xi$ where $\eta>1$ and $\xi\in\Mz$,
then $\albar=\eta$. The notion of $\tau$-localization enhanced with multiplicative decomposition turns out to be the appropriate tool 
for the purposes of the present article, whereas general $\tau$-fine-localization will re-enter the picture through the notion of closedness
in a subsequent article \cite{W}.  

\subsection{Operators related to connectivity components}
The function $\log$ ($\logend$) is described in $\Tt$-notation in Lemma 2.13 of \cite{W07c}, and for $\be=\thtt(\eta)\in\Tt$
where $\eta<\Om_1$ we have 
\begin{equation}\label{logred}
\log((1/\tau)\cdot\be)=\left\{\begin{array}{l@{\quad}l}
\eta+1&\mbox{ if }\eta=\eps+k\mbox{ where }\eps\in\Ez^{>\tau}, k<\om\\[2mm]
\eta&\mbox{ otherwise.}
\end{array}\right.
\end{equation}
The foregoing distinction reflects the property of $\tht$-functions to omit fixed points.

The operators $\iotal$ indicating the fixed-point level, $\zetal$ displaying the degree of limit point thinning, and their combination
$\latal$ measuring closure properties of ordinals $\al\in\Tt$ are as in Definitions 2.14 and 2.18 of \cite{W07c}, which also reviews the notion of 
base transformation $\pist$ and its smooth interaction with these operators.    

The operator $\latal$ already played a central role in the analysis of $\Ronepl$-patterns as it displays the number of $\le_1$-connectivity components that are $\le_1$-connected to the component with index $\al$ in a setting of relativization $\tau\in\Ezone$, 
cf.\ Lemma 2.31 part (a) of \cite{W07c}. It turns out that $\latal$ plays a similar role in $\Rtwo$, see below.

In order to avoid excessive repetition of formal definitions from \cite{CWc} we continue to describe operators and functions introduced for analysis of $\Ctwo$ in \cite{CWc} in terms of their meaning in the context of $\Ctwo$.
Those (relativized) $\leo$-components, the enumeration index of which is an epsilon number, give rise to infinite $\leo$-chains, 
along which new $\letwo$-components arise.
Omitting from these $\leo$-chains those elements that have a $\letwo$-predecessor in the chain and enumerating the remaining elements, we obtain
the so-called $\nu$-functions, see Definition 4.4 of \cite{CWc}, and the $\mu$-operator provides the length of such enumerations up to the \emph{final newly arising $\letwo$-component}, cf.\ the remark before Definition 4.4 of \cite{CWc}. This {\it terminal point} on a main line, at which the largest
newly arising $\letwo$-component originates, is crucial for understanding the structure $\Ctwo$. Note
that in general the terminal point has an infinite increasing continuation in the $\leo$-chain under consideration, leading to $\letwo$-components which have isomorphic copies below, i.e.\ which are 
{\it not} new.
Recall Convention 2.9 of \cite{W07c}.
\begin{defi}[3.4 of \cite{CWc}]
Let $\tau\in\Ezone$ and $\al\in(\tau,\tauinf)\cap\Ez$, say
$\al=\thtt(\De+\eta)$ where $\eta<\Om_1$ and $\De=\Omega_1\cdot(\la+k)$ such that $\la\in\singleton{0}\cup\Lim$ and $k<\om$.
We define
\[\mutal:=\om^{\iotal(\la)+\chial(\iotal(\la))+k}.\]\index{$\mutal$}
\end{defi}
The $\chi$-indicator occurring above is given in Definition 3.1 of \cite{CWc} and indicates whether the maximum $\le_2$-component 
starting from an ordinal on the infinite $\leo$-chain under consideration itself $\leo$-reconnects to that chain which we called a {\it main line}. 
The question remains which $\leo$-component starting from such a point on a main line is the largest that is also $\le_2$-connected to it. 
This is answered by the $\varrho$-operator: 
\begin{defi}[3.9 of \cite{CWc}]
Let $\al\in\Ez$, $\be<\al^\infty$, and $\la\in\singleton{0}\cup\Lim$, $k<\om$ be such that $\logend(\be)=\la+k$.
We define
\[\rhoalbe:=\al\cdot(\la+k\minusp\chial(\la)).\]
\end{defi}
Now, the terminal point on a main line, given as, say, $\nu^{\tauvec,\al}_{\mutal}$ with
a setting of relativization $\tauvec=(\tau_1,\ldots,\tau_n)$ that will be discussed later and $\tau=\tau_n$, $\al\in\Ez^{>\tau}$,
connects to $\latal$-many $\leo$-components. The following lemma is a direct consequence of the respective definitions.

\begin{lem}[3.12 of \cite{CWc}]\label{lamurholem}
Let $\tau\in\Ezone$ and $\al=\thtt(\De+\eta)\in(\tau,\tauinf)\cap\Ez$. Then we have
\begin{enumerate}
\item $\iotal(\De)=\rhoalmutal$ and hence $\latal=\rhoalmutal+\zetal$.
\item\label{mainlinecondpart} $\rhoalbe\le\latal$ for every $\be\le\mutal$. For $\be<\mutal$ such that
\footnote{This condition is missing in \cite{CWc}. However, that inequality was only applied under this condition, cf.\ Def.\ 5.1 and L.\ 5.7 of \cite{CWc}.}  
$\chial(\be)=0$ we even have $\rhoalbe+\al\le\latal$.
\item If $\mutal<\al$ we have $\mutal<\al\le\latal<\al^2$, while otherwise \[\max\left((\mutal+1)\cap\Ez\right)=\max\left((\latal+1)\cap\Ez\right).\]
\item If $\latal\in\Ez^{>\al}$ we have $\mutal=\latal\cdot\om$ in case of $\chial(\latal)=1$ and $\mutal=\latal$ otherwise.
\end{enumerate}
\end{lem}
Notice that we have $\mutal=\iotal(\De)=\mc(\latal)$ whenever $\mutal\in\Ez^{>\al}$.

\section{Tracking sequences and their evaluation} 

\subsection{Maximal and minimal \boldmath$\mu$\unboldmath-coverings} 

The following sets of sequences are crucial for the description of settings of relativization, which in turn is the key to
understanding the structure of connectivity components in $\Ctwo$.  

\begin{defi}[4.2 of \cite{CWc}]
Let $\tau\in\Ezone$. A nonempty sequence $(\ale,\ldots,\aln)$ of ordinals in the interval $[\tau,\tauinf)$ is called a $\tau$-tracking sequence\index{$\tau$-tracking sequence} if 
\begin{enumerate}
\item $(\ale,\ldots,\alnmin)$ is either empty or a strictly increasing sequence of epsilon numbers greater than $\tau$.
\item $\aln\in\Hz$, $\aln>1$ if $n>1$.
\item $\alie\le\mutali$ for every $i\in\singleton{1,\ldots,n-1}$.
\end{enumerate} 
By $\TSt$\index{$\tst$@$\TSt$} we denote the set of all  $\tau$-tracking sequences. Instead of $\TSe$ we also write $\TS$.\index{$\tst$@$\TSt$!$\TS$}
\end{defi}
According to Lemma 3.5 of \cite{CWc} the length of a tracking sequence is bounded in terms of the largest index of $\tht$-functions in the term representation of the first element of the sequence.

\begin{defi}\label{mucovering}
Let $\tau\in\Ezone$, $\al\in\Ez\cap(\tau,\tauinf)$, and $\be\in\Hz\cap(\al,\alinf)$.
A sequence $(\al_0,\dots,\al_{n+1})$ where $\al_0=\al$, $\al_{n+1}=\be$, 
$(\al_1,\ldots,\al_{n+1})\in\TSal$, and $\al<\al_1\le\mutal$
is called a {\it $\mu$-covering from $\al$ to $\be$}. 
\end{defi}

\begin{lem}\label{mucovloc}
Any $\mu$-covering from $\al$ to $\be$ is a subsequence of the $\al$-localization of $\be$.
\end{lem}
{\bf Proof.} Let $(\al_0,\ldots,\alne)$ be a $\mu$-covering from $\al$ to $\be$.
Stepping down from $\alne$ to $\al_0$, Lemmas 3.5 of \cite{CWc} and 4.9, 6.5 of  \cite{W07a} apply to
show that the $\al$-localization of $\be$ is the successive concatenation of the $\ali$-localization of $\alie$ for $i=1,\ldots,n$, modulo translation between the respective
settings of relativization.
\qed

\begin{defi}\label{maxminmucov}
Let $\tau\in\Ezone$.
\begin{enumerate}
\item For $\al\in\Hz\cap(\tau,\tauinf)$ we define $\maxmucovtau(\al)$ to be the longest subsequence $(\ale,\ldots,\alne)$ of the $\tau$-localization of $\al$ 
which satisfies $\tau<\ale$, $\alne=\al$, and which is {\it $\mu$-covered}, i.e.\ which satisfies $\alie\le\mutali$ for $i=1,\ldots,n$.
\item For $\al\in\Ez\cap(\tau,\tauinf)$ and $\be\in\Hz\cap(\al,\alinf)$ we denote the shortest subsequence $(\be_0,\be_1,\ldots,\be_{n+1})$
of the $\al$-localization of $\be$ which is a $\mu$-covering from $\al$ to $\be$ by $\minmucoval(\be)$, if such sequence exists. 
\end{enumerate}
\end{defi}
We recall the notion of the tracking sequence of an ordinal, for greater clarity only for 
multiplicative principals at this stage. 
\begin{defi}[cf.\ 3.13 of \cite{CWc}]\label{trsofmzdefi}
Let $\tau\in\Ezone$ and $\al\in\Mz\cap(\tau,\tauinf)$ with $\tau$-localization $\tau=\al_0,\ldots,\al_n=\al$.
The tracking sequence of $\al$ above $\tau$\index{tracking sequence}, $\trst(\al)$\index{$\trstmz$}, is defined as follows.
If there exists the largest index $i\in\{1,\ldots,n-1\}$ such that $\al\le\mutali$, then 
\[\trst(\al):=\trst(\ali)^\frown(\al),\]
otherwise $\trst(\al):=(\al)$.
\end{defi}

\begin{defi}\label{mts}
Let $\tau\in\Ezone$, $\al\in\Ez\cap(\tau,\tauinf)$, $\be\in\Hz\cap(\al,\alinf)$, and let $\al=\al_0,\ldots,\al_{n+1}=\be$ be the $\al$-localization 
of $\be$.
If there exists the least index $i\in\{0,\ldots,n\}$ such that $\ali<\be\le\mutali$, then 
\[\mtsal(\be):=\mtsal(\ali)^\frown(\be),\]
otherwise $\mtsal(\be):=(\al)$.
\end{defi}
Note that $\mtsal(\be)$ reaches $\be$ if and only if it is a $\mu$-covering from $\al$ to $\be$.

\begin{lem}\label{covcharlem} Fix $\tau\in\Ezone$.
\begin{enumerate}
\item For $\al\in\Hz\cap(\tau,\tauinf)$ let $\maxmucovtau(\al)=(\ale,\ldots,\alne)=\alvec$. If $\ale<\al$ then $\alvec$ is a $\mu$-covering from $\ale$ to $\al$
and $\mtsale(\al)\subseteq\alvec$.
\item If $\al\in\Mz\cap(\tau,\tauinf)$ then $\maxmucovtau(\al)=\trst(\al)$.
\item Let $\al\in\Ez\cap(\tau,\tauinf)$ and $\be\in\Hz\cap(\al,\alinf)$. Then $\minmucoval(\al)$ exists if and only if $\mtsal(\be)$ is a $\mu$-covering from $\al$ to $\be$,
in which case these sequences are equal, characterizing the lexicographically maximal $\mu$-covering from $\al$ to $\be$.
\end{enumerate}
\end{lem}
{\bf Proof.}  These are immediate consequences of the definitions.
\qed

\noindent Recall Definition 3.16 from \cite{CWc}, which for $\tau\in\Ezone$, $\al\in\Ez\cap(\tau,\tauinf)$ defines $\alhat$ to be the minimal $\ga\in\Mz^{>\al}$ 
such that $\trsal(\ga)=(\ga)$ and $\mutal<\ga$. 

\begin{lem}\label{mtshatlem}
Let $\tau\in\Ezone$, $\al\in\Ez\cap(\tau,\tauinf)$, and $\be\in\Mz\cap(\al,\alinf)$. Then $\mtsal(\be)$ is  a $\mu$-covering from $\al$ to $\be$
 if and only if  $\be<\alhat$. This holds if and only if for $\trsal(\be)=(\be_1,\ldots,\be_m)$ we have $\be_1\le\mutal$.
\end{lem}
{\bf Proof.} 
Suppose first that $\mtsal(\be)$ is a $\mu$-covering from $\al$ to $\be$. Then $\al$ is an element of the $\tau$-localization of $\be$, and modulo
term translation we obtain the $\tau$-localization of $\be$ by concatenating the $\tau$-localization of $\al$ with the $\al$-localization of $\be$.
By Lemma \ref{covcharlem} we  therefore have $\mtsal(\be)\subseteq(\al)^\frown\trsal(\be)$ where $\be_1\le\mual$. 
Let $\ga\in\Mz\cap(\al,\be]$ be given. Then by Lemma 3.15 of \cite{CWc} we have \[(\ga_1,\ldots,\ga_k):=\trsal(\ga)\kglex\trsal(\be),\]
so $\ga_1\le\be_1\le\mual$, and hence $\be<\alhat$.

Toward proving the converse, suppose that $\be<\alhat$. 
We have  $\trsal(\be_1)=(\be_1)$, so $\be_1\le\mual$ since $\be_1<\alhat$. This implies that $\mtsal(\be)$ reaches $\be$ as a subsequence of $(\al)^\frown\trsal(\be)$.
\qed

\begin{defi}[4.3 of \cite{CWc}] 
Let $\tau\in\Ezone$. A sequence $\alvec$ of ordinals below $\tauinf$ is a $\tau$-reference sequence\index{$\tau$-reference sequence} if 
\begin{enumerate}
\item $\alvec=()$ or
\item $\alvec=(\ale,\ldots,\aln)\in\TSt$ such that $\aln\in\Ez^{>\alnmin}$ (where $\alnod:=\tau$).
\end{enumerate}
We denote the set of $\tau$-reference sequences by $\RSt$.\index{$\rst$@$\RSt$} In case of $\tau=1$ we simply write $\RS$\index{$\rst$@$\RSt$!$\RS$} and call its elements reference sequences.\index{reference sequence}
\end{defi}

\begin{defi}[c.f.\ 4.9 of \cite{CWc}] 
For $\ga\in\Mz\cap\oneinf$ and $\eps\in\Ez\cap\oneinf$ let $\sk_\ga(\eps)$\index{$\sk$} be the maximal sequence $\de_1,\ldots,\de_l$ such that (setting $\de_0:=1$)
\begin{enumerate}
\item $\de_1=\eps$ and 
\item if $i\in\singleton{1,\ldots,l-1}\andsp\de_i\in\Ez^{>\de_{i-1}}\andsp\ga\le\mu_{\de_i}$, then 
$\de_{i+1}=\overline{\mu_{\de_i}\cdot\ga}$.
\end{enumerate}
\end{defi}

\noindent{\bf Remark(\cite{CWc}).} Lemma 3.5 of \cite{CWc} guarantees that the above definition terminates.
We have $(\de_1,\ldots,\de_{l-1})\in\RS$ and $(\de_1,\ldots,\de_l)\in\TS$. Notice that $\ga\le\de_i$ for $i=2,\ldots,l$.

\begin{defi}\label{hgamalbe}
Let $\alvec^\frown\be\in\RS$ and $\ga\in\Mz$.
\begin{enumerate}
\item If $\ga\in(\be,\behat)$, let $\mtsbe(\ga)=\etavec^\frown(\eps,\ga)$ and define
\[\hga(\alvecbe):=\alvec^\frown\etavec^\frown\skga(\eps).\]
\item If $\ga\in(1,\be]$ and $\ga\le\mube$ then 
\[\hga(\alvecbe):=\alvec^\frown\skga(\be).\]
\item If $\ga\in(1,\be]$ and $\ga>\mube$ then 
\[\hga(\alvecbe):=\alvec^\frown\be.\]
\end{enumerate}
\end{defi}

\noindent{\bf Remark.} In 1.\  let $\trsbe(\ga)=:(\ga_1,\ldots,\ga_m)$. Then we have $\ga_1\le\mube$ according to Lemma 3.17 of \cite{CWc}, 
so that $\mtsbe(\ga)$ reaches $\ga$, $\be\le\eps<\ga\le\mu_\eps$ and $\be<\mube$.
In 2.\  we have $\skga(\be)=(\be, \overline{\mube\cdot\ga})$ with $\ga\le \overline{\mube\cdot\ga}\le\be$ in case of $\mube\le\be$.
In 3.\ we have $\mube<\ga\le\be$ and hence $\skga(\be)=(\be)$.

\begin{lem}\label{hgalem}
Let $\alvec^\frown\be\in\RS$ and $\ga\in\Mz$.
Then $\hga(\alvecbe)$ is of a form $\alvec^\frown\etavec^\frown\skga(\eps)$ where
$\etavec=(\eta_1,\ldots,\eta_r)$, $r\ge 0$, $\eta_1=\be$, $\eta_{r+1}:=\eps$,
and $\skga(\eps)=(\de_1,\ldots,\de_{l+1})$, $l\ge 0$, with $\de_1=\eps$.
We have \[\lf(\de_{l+1})\ge\ga,\]
and for $\tauvecsi\in\TS$, where $\tauvec=(\tau_1,\ldots,\tau_s)$ and $\tau_{s+1}:=\si$, 
such that $\alvecbe\subseteq\tauvecsi$ and $\hga(\alvecbe)\klex\tauvecsi$
we either have 
\begin{enumerate}
\item $\tauvec=\alvec^\frown\etavec^\frown\devec_{\restriction_i}$ for some $i\in\{1,\ldots,l+1\}$ and
$\si=_\NF\de_{i+1}\cdot\sipr$ for some $\sipr<\ga$, setting $\de_{l+2}:=1$, or
\item $\tauvec_{\restriction_{s_0}}=\alvec^\frown\etavec_{\restriction_{r_0}}$ for some 
$r_0\in[1,r]$ and $s_0\le s$ such that $\eta_{r_0+1}<\tau_{s_0+1}$, in which case we have
$\mu_{\tau_j}<\ga$ for all $j\in\{s_0+1,\ldots,s\}$, and $\mu_\si<\ga$ if $\si\in\Ez^{>\tau_s}$.
\end{enumerate}
\end{lem}
{\bf Proof.}  Suppose first that $\tauvec$ is a maximal initial segment $\alvec^\frown\etavec^\frown\devec_{\restriction_i}$ for some $i\in\{1,\ldots,l+1\}$.

In the case $i=l+1$ we have $\de_{l+1}\in\Ez^{>\de_l}$ and $\mu_{\de_{l+1}}<\ga\le\de_{l+1}$, so $\si\le\mu_{\de_{l+1}}<\ga$, and we also observe that
$\tauvecsi$ could not be extended further.
Now suppose that $i\le l$. Then $\de_{i+1}<\si\le\mu_{\de_i}$, and since $\ga\le\de_{i+1}=\overline{\mu_{\de_i}\cdot\ga}\le\mu_{\de_i}$, we obtain 
$\si=\de_{i+1}\cdot\sipr$ for some $\sipr\in(1,\ga)$.

Otherwise $\tauvec$ must be of the form given in part 2 of the claim. This implies $\ga\in\Mz\cap(\be,\behat)$ and $\mtsbe(\ga)=\etavec^\frown(\eps,\ga)$.
Let us assume, toward contradiction, there existed a least $j\in\{s_0+1,\ldots,s+1\}$ such that $\tau_j\in\Ez^{>\tau_{j-1}}$ and $\ga\le\mu_{\tau_j}$.
Then ${\etavec_{\restriction_{r_0}}}^\frown(\tau_{s_0+1},\ldots,\tau_j,\ga)$ is a $\mu$-covering from $\be$ to $\ga$, hence by part 3 of Lemma \ref{covcharlem}
it must be lexicographically less than or equal to $\mtsbe(\ga)$ and therefore ${\etavec_{\restriction_{r_0}}}^\frown(\tau_{s_0+1},\ldots,\tau_j)\kglex\etavec^\frown\eps$: contradiction.
\qed

\begin{cor}\label{hgacor}
For $\alvec^\frown\be\in\RS$ and $\ga,\de\in\Mz$ such that $\de\in(1,\ga)$ we have
\[\hga(\alvecbe)\kglex\hde(\alvecbe).\]
\end{cor}

\subsection{Evaluation}

\begin{defi}\label{odef}
Let $\alvecbe\in\TS$, where $\alvec=(\ale,\ldots,\aln)$, $n\ge 0$, $\be=_\MNF\be_1\cdot\ldots\cdot\be_k$,
and set $\al_0:=1$, $\alne:=\be$, $h:=\hte(\ale)+1$, and
$\gavec_i:=\trs^{\alimin}(\ali)$, $i=1,\ldots,n$,
\[\gavec_{n+1}:=\left\{\begin{array}{ll}
            (\be)&\mbox{if }\be\le\aln\\[2mm]
            \trsaln(\be_1)^\frown\be_2&\mbox{if } k>1,\be_1\in\Ez^{>\aln}\andsp\be_2\le\mu_{\be_1}\\[2mm]
            \trsaln(\be_1)&\mbox{otherwise,}
            \end{array}\right.\]
and write $\gavec_i=(\ga_{i,1},\ldots,\ga_{i,m_i})$, $i=1,\ldots,n+1$.
Then define
\[\lSeq(\alvecbe):=(m_1,\ldots,m_{n+1})\in[h]^{\le h}.\]
Let $\bepr:=1$ if $k=1$ and $\bepr:=\be_2\cdot\ldots\cdot\be_k$ otherwise. 
We define $\ov(\alvecbe)$ recursively in $\lSeq(\alvecbe)$, as well as auxiliary parameters 
$n_0(\alvecbe)$ and $\ga(\alvecbe)$, which are set to $0$ where not defined explicitly.
\begin{enumerate}
\item $\ov((1)):=1$.
\item If $\be_1\le\aln$, then $\ov(\alvecbe):=_\NF\ov(\alvec)\cdot\be$.
\item If $\be_1\in\Ez^{>\aln}$, $k>1$, and $\be_2\le\mu_{\be_1}$, then
set $n_0(\alvecbe):=n+1$, $\ga(\alvecbe):=\be_1$, and define
\[\ov(\alvecbe):=_\NF\ov(\hop_{\be_2}(\alvec^\frown\be_1))\cdot\bepr.\] 
\item Otherwise. Then setting
\[n_0:=n_0(\alvecbe):=\max\left(\{i\in\{1,\ldots,n+1\}\mid m_i>1\}\cup\{0\}\right),\]
define 
\[\ov(\alvecbe):=_\NF\left\{\begin{array}{ll}
            \be&\mbox{if } n_0=0\\[2mm]
            \ov(\hop_{\be_1}({\alvec_{\restriction_{n_0-1}}}^\frown\ga))\cdot\be&\mbox{if } n_0>0,
            \end{array}\right.\] 
where $\ga:=\ga(\alvecbe):=\ga_{n_0,m_{n_0}-1}$.
\end{enumerate}
\end{defi}

\noindent{\bf Remark.}  As indicated in writing $=_\NF$ in the above definition, we obtain terms in multiplicative normal form
denoting the values of $\ov$. The {\bf fixed points of $\ov$}, i.e.\ those $\alvecbe$ that satisfy $\ov(\alvecbe)=\be$ are therefore
characterized by 1.\ and 4.\ for $n_0=0$.

Recall Definition 3.13 of \cite{CWc}, extending Definition \ref{trsofmzdefi} to additive principal numbers that are not multiplicative
principal ones.

\begin{defi}[cf.\ 3.13 of \cite{CWc}]\label{trsofhzdefi}
Let $\tau\in\Ezone$ and $\al\in[\tau,\tauinf)\cap\Hz$.
The tracking sequence of $\al$ above $\tau$\index{tracking sequence}, $\trst(\al)$\index{$\trsthz$}, is defined as in Definition \ref{trsofmzdefi} 
if $\al\in\Mz^{>\tau}$, and otherwise recursively in the multiplicative decomposition of $\al$ as follows.
\begin{enumerate}
\item If $\al\le\tau^\om$ then $\trst(\al):=(\al)$.
\item Otherwise. Then $\albar\in[\tau,\al)$ and $\al=_\mNF\albar\cdot\be$ for some $\be\in\Mz^{>1}$.
      Let $\trst(\albar)=(\ale,\ldots,\aln)$ and set $\al_0:=\tau$.\footnote{As verified in part 2 of the lemma below we have 
      $\be\le\al_n$.}
      \begin{enumerate}
      \item[2.1.] If $\aln\in\Ez^{>\al_{n-1}}$ and $\be\le\mutaln$ then $\trst(\al):=(\ale,\ldots,\aln,\be)$.
      \item[2.2.] Otherwise. For $i\in\singleton{1,\ldots,n}$ let $(\be^i_1,\ldots,\be^i_{m_i})$ be $\trs^{\ali}(\be)$
            provided $\be>\ali$, and set $m_i:=1$, $\be^i_1:=\ali\cdot\be$ if $\be\le\ali$.
            We first define the critical index
            \[i_0(\al)=i_0:=\max\left(\singleton{1}\cup\set{j\in\singleton{2,\ldots,n}}{\be^j_1\le\mu^\tau_{\al_{j-1}}}\right).\]
            Then $\trst(\al):=(\ale,\ldots,\al_{i_0-1},\be^{i_0}_1,\ldots,\be^{i_0}_{m_{i_0}})$.
      \end{enumerate}
\end{enumerate}            
Instead of $\trs^1(\al)$ we also simply write $\trs(\al)$.\index{$\trst$!$\trs$}
\end{defi}

\begin{lem}\label{trsestimlem} 
If in the above definition, part 2.2, we have $\be>\alinod$, then for all $j\in(i_0,\ldots,n]$ we have \[\be^{i_0}_1\le\alj\le\mu_{\al_{j-1}},\] 
in particular $\be^{i_0}_1\le\mu_\alinod$.
\end{lem}
{\bf Proof.} Assume toward contradiction that there exists the maximal $j\in(i_0,\ldots,n]$ such that $\alj<\be^{i_0}_1$.
Since $\be^{i_0}_1\le\be\le\aln$ we have $j<n$ and obtain \[\alinod<\alj<\be^{i_0}_1\le\al_{j+1}\le\mu_\alj,\]
implying that $\alj\in\trs^\alinod(\be)$, contradicting the minimality of $\be^{i_0}_1$ in $\trs^\alinod(\be)$. \qed

\begin{lem}[3.14 of \cite{CWc}]
Let $\tau\in\Ezone$ and $\al\in[\tau,\tauinf)\cap\Hz$. Let further $(\ale,\ldots,\aln)$ be
$\trst(\al)$, the tracking sequence of $\al$ above $\tau$.
\begin{enumerate}
\item If $\al\in\Mz$ then $\aln=\al$ and $\trst(\ali)=(\ale,\ldots,\ali)$ for $i=1,\ldots,n$.
\item If $\al=_\mNF\eta\cdot\xi\not\in\Mz$ then $\aln\in\Hz\cap[\xi,\al]$ and $\aln=_\mNF\alnbar\cdot\xi$.
\item $(\ale,\ldots,\al_{n-1})$ is either empty or a strictly increasing sequence of epsilon numbers in the interval $(\tau,\al)$. 
\item For $1\le i\le n-1$ we have $\alie\le\mutali$, and if $\ali<\alie$ then $(\ale,\ldots,\alie)$ is a subsequence of the $\tau$-localization of $\alie$.
\end{enumerate}
\end{lem}
{\bf Proof.} The proof proceeds by straightforward induction along the definition of $\trst(\al)$, i.e. along the length of the 
$\tau$-localization of multiplicative principal numbers and
the number of factors in the multiplicative decomposition of additive principal numbers. In part 4 Lemma 6.5 of \cite{W07a} and the previous remark apply.\qed

\begin{lem}[3.15 of \cite{CWc}]\label{citedinjtrslem}
Let $\tau\in\Ezone$ and $\al,\ga\in[\tau,\tauinf)\cap\Hz$, $\al<\ga$. Then we have 
\[\trst(\al)\klex\trst(\ga).\] 
\end{lem}
{\bf Proof.} The proof given in \cite{CWc} is in fact an induction along the inductive definition of $\trst(\ga)$ with a subsidiary induction along the 
inductive definition of $\trst(\al)$.
\qed

\begin{theo}\label{thma}
For all $\al\in\Hz\cap\oneinf$ we have \[\ov(\trs(\al))=\al.\] 
\end{theo}
{\bf Proof.}  The theorem is proved by induction along the inductive definition of $\trs(\al)$.\\[2mm]
{\bf Case 1:} $\al\in\Mz$. Then $\lSeq(\trs(\al))=(1,\ldots,1)$ and hence $\ov(\trs(\al))=\al$ immediately by definition.\\[2mm]
{\bf Case 2:} $\al=_\NF\albar\cdot\be\in\Hz-\Mz$. Let $\trs(\albar)=:(\ale,\ldots,\aln)$ and $\al_0:=1$. By the i.h.\ $\ov(\alvec)=\albar$.
We have $\be\le\lf(\aln)\le\aln$, $n\ge 1$.\\[2mm]
{\bf Subcase 2.1:} $\aln\in\Ez^{>\alnmin}\andsp\be\le\mualn$. Then $\trs(\al)=\alvec^\frown\be$, and since $\be\in\Mz^{\le\aln}$ according to the definition of $\ov$ 
we obtain $\ov(\alvecbe)=\ov(\alvec)\cdot\be=\albar\cdot\be=\al$.\\[2mm]
{\bf Subcase 2.2:} Otherwise. Let $(\be^i_1,\ldots,\be^i_{m_i})$  for $i=1,\ldots,n$ as well as the index $i_0$ be defined as in case 2.2 of Definition \ref{trsofhzdefi},
so that $\trs(\al)=(\ale,\ldots,\al_{i_0-1},\be^{i_0}_1,\ldots,\be^{i_0}_{m_{i_0}})$.\\[2mm]
{\bf 2.2.1:} $i_0=n$. Then we have $\trs(\al)=(\ale,\ldots,\alnmin,\aln\cdot\be)$, and using the i.h.\ we obtain $\ov(\trs(\al))=\ov(\alvec)\cdot\be=\al$.\\[2mm]
{\bf 2.2.2:} $i_0<n$ and $\be\le\alinod$. Then we have $\alinod\in\Ez^{>\al_{i_0-1}}$, $\alinod\cdot\be\le\mu_{\al_{i_0-1}}$, and 
$\trs(\al)=(\ale,\ldots,\al_{i_0-1},\alinod\cdot\be)$. It follows that for all $j\in(i_0,n]$ we have $\be\le\alj$ and $\alj\cdot\be>\mu_{\al_{j-1}}$, hence $\be\le\mu_\alinod$
and thus \[\ov(\trs(\al))=\ov(\hbe(\alvec_{\restriction_{i_0}}))\cdot\be.\] The sequence $\hbe(\alvec_{\restriction_{i_0}})$ is of the form 
${\alvec_{\restriction_{i_0-1}}}^\frown\devec$ where $\devec:=\skbe(\alinod)$. Since $\mu_\alinod<\al_{i_0+1}\cdot\be$ we have 
$\overline{\mu_\alinod\cdot\be}=\al_{i_0+1}=\de_2$. In the case $i_0+1=n$ we obtain $\devec=(\alnmin,\aln)$, hence $\hbe(\alvec_{\restriction_{i_0}})=\alvec$,
otherwise we iterate the above argumentation to see that $\devec=(\alinod,\ldots,\aln)$. Hence $\ov(\trs(\al))=\ov(\alvec)\cdot\be$ as desired.\\[2mm]
{\bf 2.2.3:} $i_0<n$ and $\alinod<\be$. Then we have $\bevec^{i_0}=\trs^\alinod(\be)$ with $\be^{i_0}_1\le\mu_{\al_{i_0-1}}$ if $i_0>1$.
By Lemma \ref{trsestimlem} $\alinod$ is the immediate predecessor of $\be^{i_0}_1$ in $\trs^{\al_{i_0-1}}(\be)$.
By definition of $\ov$ we have $\ov(\trs(\al))=\ov(\hbe(\alvec_{\restriction_{i_0}}))\cdot\be$ and therefore have to show that $\hbe(\alvec_{\restriction_{i_0}})=\alvec$.
We define \[j_0:=\min\{j\in\{i_0,\ldots,n-1\}\mid\be\le\mu_\alj\},\]
which exists, because $\be\le\aln\le\mu_\alnmin$.\\[2mm]
{\bf Claim:} $\mts^\alinod(\be)=(\alinod,\ldots,\aljnod,\be)$.\\
{\bf Proof.} For every $j\in\{i_0,\ldots,j_0-1\}$ the minimality of $j_0$ implies $\be>\mualj\ge\alje$, and thus by the maximality of $i_0$ also $\be^{j+1}_1>\mualj$.
Moreover, we have \[\be^{j+1}_1\le\al_{j+2}:\] Assume otherwise and let $j$ be maximal in $\{i_0,\ldots,j_0-1\}$ such that $\be^{j+1}_1>\al_{j+2}$. Since
$\be^{j+1}_1\le\be\le\aln$ we must have $j\le n-3$. But then $\al_{j+1}<\al_{j+2}<\be^{j+1}_1\le\al_{j+3}\le\mu_{j+2}$ and hence $\al_{j+2}\in\trs^{\al_{j+1}}(\be)$,
contradicting the minimality of $\be^{j+1}_1$ in $\trs^{\al_{j+1}}(\be)$.
Therefore \[\mualj<\be^{j+1}_1\le\al_{j+2}\le\mu_{\alje},\]
which concludes the proof of the claim.\qed

It remains to be shown that $\skbe(\aljnod)=(\aljnod,\ldots,\aln)$, i.e.\ to successively check for $j=j_0,\ldots,n-1$ that $\be\le\alje\le\mualj$ and $\alje\cdot\be>\mualj$,
whence $\overline{\mualj\cdot\be}=\alje$. This concludes the verification of $\hbe(\alvec_{\restriction_{i_0}})=\alvec$ and consequently the proof of the theorem.
\qed

\begin{theo}\label{thmb}
For all $\alvecbe\in\TS$ we have \[\trs(\ov(\alvecbe))=\alvecbe.\]
\end{theo}
{\bf Proof.} The theorem is proved by induction on $\lSeq(\alvecbe)$ along the ordering $(\lSeq,\klex)$.
Let $\be=_\NF\be_1\cdot\ldots\cdot\be_k$ and set $n_0:=n_0(\alvecbe)$, $\ga:=\ga(\alvecbe)$ according to Definition \ref{odef}, which 
provides us with an $\NF$-representation of $\ov(\alvecbe)$, where in the interesting cases the i.h.\ applies to the term $\trs\left(\overline{\ov(\alvecbe)}\right)$.   
\\[2mm]
{\bf Case 1:} $n=0$ and $\be=1$. Trivial.\\[2mm]
{\bf Case 2:} $1<\be_1\le\aln$. Then $\ov(\alvecbe)=_\NF\ov(\alvec)\cdot\be$, and it is straightforward to verify the claim from the i.h.\ applied to $\alvec$ 
by inspecting case 2.1 of Definition \ref{trsofhzdefi}.\\[2mm]
{\bf Case 3:} $k>1$ with $\be_1\in\Ez^{>\aln}$ and $\be_2\le\mu_{\be_1}$. Then by definition $\ov(\alvecbe)=_\NF\ov(\hop_{\be_2}(\alvec^\frown\be_1))\cdot\bepr$
where $\bepr=(1/\be_1)\cdot\be$. According to part 2 of Definition \ref{hgamalbe} we have \[\hop_{\be_2}(\alvec^\frown\be_1)=\alvec^\frown\devec,\] where
$\devec=(\de_1,\ldots,\de_{l+1}):=\sk_{\be_2}(\be_1)$. Assume first that $k=2$. Then the maximality of the length of $\devec$ excludes the possibility
$\de_{l+1}\in\Ez^{>\de_l}\andsp\be_2\le\mu_{\de_{l+1}}$. We have $\be_2\le\be_1$ and $\be\le\mualn$. For $j\in\{2,\ldots,l+1\}$ we have 
$\be_2\le\de_j=\overline{\mu_{\de_{j-1}}\cdot\be_2}\le\mu_{\de_{j-1}}$ and $\de_j\cdot\be_2=\mu_{\de_{j-1}}\cdot\be_2>\mu_{\de_{j-1}}$.
This implies that $\trs(\ov(\alvecbe))=\alvec^\frown\be$. The claim now follows easily for $k>2$ since $\be\le\mualn$.\\[2mm]
{\bf Case 4:} Otherwise.\\[2mm]
{\bf Subcase 4.1:} $n_0=0$. Then we have $\ov(\alvecbe)=\be$, $\trs^\alimin(\ali)=(\ali)$ for $i=1,\ldots,n$, and $\trs^\aln(\be_1)=(\be_1)$, whence
$\trs(\be)=(\be)$.\\[2mm]
{\bf Subcase 4.2:} $n_0>0$. Using the abbreviation $\alvecpr:=\alvec_{\restriction{n_0-1}}$ we then have $\ov(\alvecbe)=_\NF\ov(\hop_{\be_1}(\alvecpr^\frown\ga))\cdot\be$. 
Setting $\mts^\ga(\be_1)=:(\ga_1,\ldots,\ga_{m+1})$ where $\ga_1=\ga$, $\ga_{m+1}=\be_1$, and $\sk_{\be_1}(\ga_m)=:\devec=(\de_1,\ldots,\de_{l+1})$ where $\de_1=\ga_m$,
we have \[\hop_{\be_1}(\alvecpr^\frown\ga)={\alvecpr}^\frown{\gavec_{\restriction_{m-1}}}^\frown\devec,\]
which by the i.h.\ is the tracking sequence of  $\overline{\ov(\alvec^\frown\be_1)}$. 
Assuming first that $k=1$, we now verify that the tracking sequence of $\ov(\alvecbe)$ actually is $\alvecbe$, by checking
that case 2.2 of Definition \ref{trsofhzdefi} applies, with $n_0$ playing the role of the critical index $i_0(\ov(\alvecbe))$.
Note first that the maximality of the length of $\devec$ rules out the possibility $\de_{l+1}\in\Ez^{>\de_l}\andsp\be_1\le\mu_{\de_{l+1}}$ and hence case 2.1 of 
Definition \ref{trsofhzdefi}.
According to the choice of $n_0$ and part 2 of Lemma \ref{covcharlem} we have $\trs^\ga(\be_1)=\maxmucov^\ga(\be_1)=(\al_{n_0},\ldots,\aln,\be_1)$ and
of course $\al_{n_0}\le\mu_{\al_{n_0-1}}$. Thus $n_0$ qualifies for the critical index, once we show its maximality:
Firstly, for any $i\in\{2,\ldots,m\}$ we have $\ga_i<\be_1$, and setting $\bevec^i:=\trs^{\ga_i}(\be_1)$ the assumption $\be^i_1\le\mu_{\ga_{i-1}}$ would imply
that ${\gavec_{\restriction_{i-1}}}^\frown\bevec^i$ is a $\mu$-covering from $\ga$ to $\be_1$ such that
\[\mts^\ga(\be_1)\klex{\gavec_{\restriction_{i-1}}}^\frown\bevec^i,\]
contradicting part 3 of Lemma \ref{covcharlem}. 
Secondly, for any $j\in\{2,\ldots,l+1\}$ we have $\de_j=\overline{\mu_{\de_{j-1}}\cdot\be_1}$, so $\de_j\cdot\be_1=\mu_{\de_{j-1}}\cdot\be_1>\mu_{\de_{j-1}}$.
These considerations entail \[\trs(\alvec^\frown\be_1)=\alvecpr^\frown\trs^\ga(\be_1),\]
and it is easy now to verify the claim for arbitrary $k$, again since $\be\le\mualn$.  
\qed

\begin{cor}\label{ocontcor}
$\ov$ is strictly increasing with respect to the lexicographic ordering on $\TS$ and 
continuous in the last vector component. 
\end{cor}
{\bf Proof.} The first statement is immediate from Lemma \ref{citedinjtrslem} and Theorems \ref{thma} and \ref{thmb}.
In order to verify continuity, let $\alvec=(\ale,\ldots,\aln)\in\RS$ and $\be\in\Lz^{\le\mualn}$ be given.
For any $\ga\in\Hz\cap\be$, we have $\gati:=\ov(\alvecga)<\ov(\alvecbe)=:\beti$ and 
$\alvecga=\trs(\gati)\klex\trs(\beti)=\alvecbe$. For given $\deti\in\Hz\cap(\gati,\beti)$ set
$\devec:=\trs(\deti)$, so that \[\alvecga\klex\devec\klex\alvecbe,\]
whence $\alvec\subseteq\devec$ is an initial segment. Writing $\devec=\alvec^\frown(\ze_1,\ldots,\ze_m)$
we obtain $\ga<\ze_1<\be$. For $\devecpr:=\alvec^\frown\ze_1\cdot\om$ we then have $\deti<\ov(\devecpr)<\beti$.
\qed

\noindent{\bf Remark.} Theorems \ref{thma} and \ref{thmb} establish Lemma 4.10 of \cite{CWc} in a weak theory, 
adjusted to our redefinition of $\ov$. Its equivalence with the definition in \cite{CWc} follows, since
the definition of $\trs$ has not been modified. In the next section we will continue in this way in order to obtain suitable
redefinitions of the $\kappa$- and $\nu$-functions. 

The following lemma will not be required in the sequel, but it
has been included to further illuminate the approach.
 
\begin{lem}\label{olem}
Let $\alvecbe\in\RS$, $\ga\in\Mz\cap(1,\behat)$, and let $\tauvecsi\in\TS$ be such that $\alvecbe\subseteq\tauvecsi$ and $\hga(\alvecbe)\klex\tauvecsi$.
Then we have \[\ov(\tauvecsi)=_\NF\ov(\hga(\alvecbe))\cdot\de\]
for some $\de\in\Hz\cap(1,\ga)$.
\end{lem}
{\bf Proof.} The lemma is proved by induction on $\lSeq(\tauvecsi)$, using Lemma \ref{hgalem}. 
In order to fix some notation, set $\tauvec=(\tau_1,\ldots,\tau_s)$ and $\tau_{s+1}:=\si$.
Write $\hga(\alvecbe)=\alvec^\frown\etavec^\frown\devec$ where $\etavec=(\eta_1,\ldots,\eta_r)$ and $\devec=(\de_1,\ldots,\de_{l+1})=\skga(\eps)$, 
$\de_1=\eps=:\eta_{r+1}$, and $\de_{l+2}:=1$.
Note that since $\alvecbe\in\RS$, we have $\eps\in\Ez^{>\eta_r}$. According to Lemma \ref{hgalem} either one of the two following cases applies
to $\tauvecsi$. \\[2mm]
{\bf Case 1:} $\tauvec=\alvec^\frown\etavec^\frown\devec_{\restriction_i}$ for some $i\in\{1,\ldots,l+1\}$ and $\si=_\NF\de_{i+1}\cdot\ze$ for some 
$\ze\in\Hz\cap(1,\ga)$. Let $\ze=_\MNF\ze_1\cdot\ldots\cdot\ze_j$.\\[2mm]
{\bf Subcase 1.1:} $i=l+1$. Then $\hga(\alvecbe)=\tauvec$, $\de_{l+1}\in\Ez^{>\de_l}$, and 
$\si\le\mu_{\de_{l+1}}<\ga\le\de_{l+1}$. According to the definition of $\ov$ we have 
$\ov(\tauvecsi)=\ov(\hga(\alvecbe))\cdot\si$.\\[2mm]
{\bf Subcase 1.2:} $i\le l$ where $\de_{i+1}\in\Ez^{>\de_i}$ and $\ze_1\le\mu_{\de_{i+1}}$.
By definition, $\ov(\tauvecsi)=\ov(\hop_{\ze_1}(\tauvec^\frown\de_{i+1}))\cdot\ze$.
As $\tauvec^\frown\de_{i+1}$ is an initial segment of $\hga(\alvecbe)$ and $\ze_1<\ga$, we have 
\[\hga(\alvecbe)=\hga(\tauvec^\frown\de_{i+1})\kglex\hop_{\ze_1}(\tauvec^\frown\de_{i+1})\]
by Corollary \ref{hgacor}. The claim now follows by the i.h.\ (if necessary) applied 
to $\hop_{\ze_1}(\tauvec^\frown\de_{i+1})$.\\[2mm]
{\bf Subcase 1.3:} Otherwise. Here we must have $i=l$, since if $i<l$ it follows that $\de_{i+1}\in\Ez^{>\de_i}$
and $\ze_1<\ga\le\overline{\mu_{\de_{i+1}}\cdot\ga}=\de_{i+2}\le\mu_{\de_{i+1}}$, which has been covered by the
previous subcase. We therefore have $\ov(\tauvecsi)=\ov(\hga(\alvecbe))\cdot\ze$.\\[2mm]
{\bf Case 2:} $\tauvec_{\restriction_{s_0}}=\alvec^\frown\etavec_{\restriction_{r_0}}$ for some $r_0\in[1,r]$, $s_0\le s$, $\eta_{r_0+1}<\tau_{s_0+1}$,
whence according to Lemma \ref{hgalem} we have $\si\le\mu_{\tau_s}<\ga$, $\mu_{\tau_j}<\ga$ for $j=s_0+1,\ldots,s$, and $\mu_\si<\ga$ if $\si\in\Ez^{>\tau_s}$.
Let $\si=_\MNF\si_1\cdot\ldots\cdot\si_k$ and $\sipr:=(1/\si_1)\cdot\si$.

In the case $s_0<s$ we have $\hga(\alvecbe)\klex\tauvec$, and noting that $\si<\ga$, the i.h.\ straightforwardly applies to $\tauvec$. Let us therefore assume that $s_0=s$,
whence $\tauvec=\alvec^\frown\etavec_{\restriction_{r_0}}$ and $\eta_{r_0+1}<\si$. Then we have $\etavec^\frown\eps\klex\tauvecsi$ and 
$\eta_{r_0}<\eta_{r_0+1}<\si\le\mu_{\eta_{r_0}}<\ga$.\\[2mm]
{\bf Subcase 2.1:} $k>1$ where $\si_1\in\Ez^{>\eta_{r_0}}$ and $\si_2\le\mu_{\si_1}$.
Then $\ov(\tauvecsi)=\ov(\hop_{\si_2}(\tauvec^\frown\si_1))\cdot\sipr$ and 
$\ga>\si_1\ge\eta_{r_0+1}\in\Ez^{>\eta_{r_0}}$. In the case $\si_1>\eta_{r_0+1}$ the i.h.\ applies to
$\hop_{\si_2}(\tauvec^\frown\si_1)$. Now assume $\si_1=\eta_{r_0+1}$. As in Subcase 1.2 we obtain
\[\hga(\alvecbe)=\hga(\tauvec^\frown\si_1)\kglex\hop_{\si_2}(\tauvec^\frown\si_1),\]
and (if necessary) the i.h.\ applies to $\hop_{\si_2}(\tauvec^\frown\si_1)$.\\[2mm]
{\bf Subcase 2.2:} Otherwise, i.e.\ Case 4 of Definition \ref{odef} applies. Let $n_0:=n_0(\tauvecsi)$.
We first assume that $k=1$, i.e.\ $\si\in\Mz\cap(\eta_{r_0+1},\ga)$.\\[2mm]
{\bf 2.2.1:} $n_0\le s$. We obtain $\ov(\tauvec)=\xi\cdot\tau_s$ and $\ov(\tauvecsi)=\xi\cdot\si$ where $\xi=1$ if
$n_0=0$ and $\xi=\hop_{\tau_s}({\tauvec_{\restriction_{n_0-1}}}^\frown\ga(\tauvec))$.
The $<,\klex$-order isomorphism between $\TS$ and $\Hz^{<\oneinf}$ established by
Lemma \ref{citedinjtrslem} and Theorems \ref{thma} and \ref{thmb} yields
\[\ov(\tauvec)<\ov(\hga(\alvecbe))<\ov(\tauvecsi),\]
and hence the claim.\\[2mm]
{\bf 2.2.2:} $n_0=s+1$. We have $\si\in\Mz\cap(\eta_{r_0+1},\ga)$. Let $\ze$ be the immediate predecessor
of $\si$ in $\trs^{\eta_{r_0}}(\si)$. Then $\ov(\tauvecsi)=\ov(\hop_\si(\tauvec^\frown\ze))\cdot\si$.\\[2mm]
{\bf 2.2.2.1:} $\eta_{r_0+1}<\ze$. Then the i.h.\ applies to $\tauvec^\frown\ze$.\\[2mm]
{\bf 2.2.2.2:} $\eta_{r_0+1}=\ze$. Then we argue as before, since 
\[\hga(\alvecbe)=\hga(\tauvec^\frown\eta_{r_0+1})\kglex\hop_\si(\tauvec^\frown\eta_{r_0+1}),\]
so that the i.h.\ (if necessary) applies to $\hop_\si(\tauvec^\frown\eta_{r_0+1})$.\\[2mm]
{\bf 2.2.2.3:} $\eta_{r_0+1}>\ze$. Then $\ze$ is an element of $\trs^{\eta_{r_0}}(\eta_{r_0+1})$, and
by a monotonicity argument as in 2.2.1 we obtain the claim as a consequence of 
\[\ov(\hop_\si(\tauvec^\frown\ze))<\ov(\hga(\alvecbe))<\ov(\tauvecsi).\]
This concludes the proof for $k=1$, and for $k>1$ the claim now follows easily. 
\qed

\section{Enumerating relativized connectivity components}\label{conncompsec}
Recall Definition 4.4 of \cite{CWc}. We are now going to characterize the functions $\kappa$ and $\nu$
by giving an alternative definition which is considerably less intertwined. The first step is to define
the restrictions of $\kval$ and $\nuval$ to additive principal indices. Recall part 3 of Lemma \ref{lamurholem}.

\begin{defi}\label{kappanuprincipals}
Let $\alvec\in\RS$ where $\alvec=(\ale,\ldots,\aln)$, $n\ge 0$, $\al_0:=1$.
We define $\kvalbe$ and $\nuvalbe$ for additive principal $\be$ as follows, writing $\kappa_\be$ instead of $\kappa^{()}_\be$.\\[2mm]
{\bf Case 1:} $n=0$. For $\be<\oneinf$ define \[\kappa_\be:=\ov((\be)).\]
{\bf Case 2:} $n>0$. For $\be\le \mualn$, i.e.\ $\alvec^\frown\be\in\TS$, define \[\nu^\alvec_\be:=\ov(\alvecbe).\]
$\kvalbe$ for $\be\le \laaln$ is defined by cases. If $\be\le\aln$ let $i\in\{0,\ldots,n-1\}$ be maximal such that $\ali<\be$.
If $\be>\aln$ let $\be=_\MNF\be_1\cdot\ldots\cdot\be_k$ and set $\bepr:=(1/\be_1)\cdot\be$. 
\[\kvalbe:=\left\{\begin{array}{ll}
            \ka^{{\alvec_{\restriction_i}}}_\be&\mbox{if } \be\le\aln\\[2mm]
            \ov(\alvec)\cdot\bepr&\mbox{if } \be_1=\aln\andsp k>1\\[2mm]
            \ov(\alvecbe)&\mbox{if } \be_1>\aln.
            \end{array}\right.\]
\end{defi}

\noindent{\bf Remark.}  Note that in the case $n>0$ we have the following inequalities between $\kvalbe$ and $\nuvalbe$, which are consequences
of the monotonicity of $\ov$ proved in Theorems \ref{thma} and \ref{thmb}. 
\begin{enumerate}
\item If $\be\le\aln$ then $\kvalbe\le\kvalnminaln=\ov(\alvec)$. Later we will define $\nuval_0:=\ov(\alvec)$.
\item If $\be_1=\aln$ and $k>1$ then $\kvalbe=\ov(\alvec)\cdot\bepr=\nuval_\bepr$, which is less than $\nuvalbe$ if $\be\le\mualn$.
\item Otherwise we have $\kvalbe=\nuvalbe$.
\end{enumerate}

\begin{cor}\label{kappanucor}
$\mbox{ }$
\begin{enumerate}
\item $\kappa$ and $\nu$ are strictly increasing with respect to their $\klex$-ordered arguments $\alvecbe\in\TS$.
\item Each branch $\kval$ (where $\alvec\in\RS$) and $\nuval$ (where $\alvec\in\RS-\{()\}$) is continuous at arguments $\be\in\Lz$.
\end{enumerate}
\end{cor}
{\bf Proof.} This is a consequence of Corollary \ref{ocontcor}.
\qed

We now prepare for the conservative extension of $\kappa$ and $\nu$ to their entire domain as well as
the definition of $\dpf$ which is in accordance with Definition 4.4 of \cite{CWc}.

\begin{defi}\label{Ttauvec}
Let $\tauvec\in\RS$, $\tauvec=(\tau_1,\ldots,\tau_n)$, $n\ge0$, $\tau_0:=1$.
The term system $\Ttvec$ is obtained from $\Ttn$ by successive substitution of parameters in $(\taui,\tauie)$ by their $\Tti$-representations,
for $i=n-1,\ldots,1$. The parameters $\taui$ are represented by the terms $\thtti(0)$.
The length $\ltvec(\al)$ of a $\Ttvec$-term $\al$ is defined inductively by 
\begin{enumerate}
\item $\ltvec(0):=0$,
\item $\ltvec(\be):=\ltvec(\ga)+\ltvec(\de)$ if $\be=_\NF\ga+\de$, and 
\item $\ltvec(\tht(\eta)):=\left\{\begin{array}{l@{\quad}l}
1&\mbox{ if }\quad\eta=0\\
\ltvec(\eta)+4&\mbox{ if }\quad\eta>0
\end{array}\right.$\\[2mm] 
where $\tht\in\{\tht^{\tau_i}\mid 0\le i\le n\}\cup\{\tht_{i+1}\mid i\in\N\}$.
\end{enumerate}
\end{defi}

\noindent{\bf Remark.}  Recall Equation (\ref{logred}) as well as Lemma 2.13 and Definitions 2.14 and 2.18 of \cite{W07c}.
\begin{enumerate}
\item For $\be=\tht^{\tau_n}(\De+\eta)\in\Ez$ such that $\be\le\mu_{\tau_n}$ we have
\begin{equation}\label{iotalen} \ltvec(\De)=\ltvecbe(\iota_{\tau_n,\be}(\De))<\ltvec(\be).
\end{equation} 
\item For $\be\in\Ttvec\cap\Hz^{>1}\cap\Om_1$ let $\tau\in\{\tau_0,\ldots,\tau_n\}$ be maximal such that $\tau<\be$. 
Clearly,
\begin{equation}\label{barlen} \ltvec(\bebar) < \ltvec(\be),
\end{equation}
cf.\ Subsection \ref{reflocsubsec}, and
\begin{equation}\label{zelen}
\ltvec(\zetbe) < \ltvec(\be),
\end{equation}
In case of $\be\not\in\Ez$ we have
\begin{equation}\label{loglen} \ltvec(\log(\be)), \ltvec(\log((1/\tau)\cdot\be)) < \ltvec(\be),
\end{equation}
and for $\be\in\Ez$ we have
\begin{equation}\label{lalen} 
\ltvecbe(\latbe) < \ltvec(\be).
\end{equation}
\end{enumerate}
 
Finally, the definition of the enumeration functions of relativized connectivity components can be
completed. This is easily seen to be a sound, elementary recursive definition. 

\begin{defi}[cf.\ 4.4 of \cite{CWc}]
Let $\alvec\in\RS$ where $\alvec=(\al_1,\ldots,\al_n)$, $n\ge0$, and set $\al_0:=1$. 
We define the functions
\[\kval, \dpval: \domkval\to\oneinf,\] 
\index{$\kval$}\index{$\dpval$}

\noindent where $\domkval:=\oneinf$ if $n=0$ and $\domkval:=[0,\laaln]$ if $n>0$, simultaneously by recursion on $\lalvec(\be)$, extending Definition \ref{kappanuprincipals}.
The clauses extending the definition of $\kval$ are as follows.  
\begin{enumerate}
\item $\kval_0:=0$, $\kval_1:=1$,
\item\label{kappapl} $\kvalbe:=\kvalga+\dpval(\ga)+\kvalde$ for $\be=_\NF\ga+\de$.
\end{enumerate}

\noindent $\dpval$ is defined as follows, using $\nu$ as already defined on $\TS$.
\begin{enumerate}
\item $\dpval(0):=0$, $\dpval(1):=0$, and $\dpval(\aln):=0$ in case of $n>0$,
\item $\dpval(\be):=\dpval(\de)$ if $\be=_\NF\ga+\de$,
\item\label{dpred} $\dpval(\be):=\dpf_{\alvecrestrnmin}(\be)$ if $n>0$ for $\be\in\Hz\cap(1,\aln)$,
\item for $\be\in\Hz^{>\aln}-\Ez$ let $\ga:=(1/\aln)\cdot\be$ and $\log(\ga)=_\ANF\ga_1+\ldots+\ga_m$ and set 
      \[\dpval(\be):=\kval_{\ga_1}+\dpval(\ga_1)+\ldots+\kval_{\ga_m}+\dpval(\ga_m),\]
\item\label{dpeps} for $\be\in\Ez^{>\aln}$ let $\gavec:=(\ale,\ldots,\aln,\be)$, and set        
      \[\dpval(\be):=\nuvga_{\mu^\aln_\be}+\kvga_{\laalnbe}+\dpvga(\laalnbe).\]
\end{enumerate}
\end{defi}

\begin{defi}[cf.\ 4.4 of \cite{CWc}]
Let $\alvec\in\RS$ where $\alvec=(\al_1,\ldots,\al_n)$, $n>0$, and set $\al_0:=1$. 
We define
\[\nuval:\domnuval\to\oneinf\]\index{$\nuval$}

\noindent where $\domnuval:=[0,\mu_\aln]$, 
extending Definition \ref{kappanuprincipals} and setting $\al:=\ov(\alvec)$, by
\begin{enumerate}
\item $\nuval_0:=\al$, 
\item\label{nuple} $\nuval_{\be}:=\nuval_\ga+\kval_{\rhoalnga}+\dpval(\rhoalnga)+\chialncheck(\ga)\cdot\al$ if $\be=\ga+1$,
\item $\nuval_\be:=\nuval_\ga+\kval_{\rhoalnga}+\dpval(\rhoalnga)+\nuval_{\de}$ if $\be=_\NF\ga+\de\in\Lim$.
\end{enumerate}
\end{defi}


In the sequel we want to establish the results of Lemma 4.5 of \cite{CWc} for the new definitions within a weak theory, avoiding the long transfinite induction 
used in the corresponding proof in \cite{CWc}.
Then the agreement of the definitions of $\kappa,\dpf$ and $\nu$ in \cite{CWc} and here can be shown in a weak theory as well. This includes also Lemma 4.7 of \cite{CWc}
and extends to the relativization of tracking sequences to contexts as lined out in Definition 4.13 through Lemma 4.17 of \cite{CWc}.


\begin{lem}\label{kdpmainlem}
Let $\alvec=(\ale,\ldots,\aln)\in\RS$ and set $\al_0:=1$.
\begin{enumerate}
\item Let $\ga\in\domkval\cap\Hz$. If $\ga=_\MNF\ga_1\cdot\ldots\cdot\ga_k\ge\aln$, setting
$\gapr:=(1/\ga_1)\cdot\ga$, we have 
\[(\kvalga+\dpval(\ga))\cdot\om=\left\{\begin{array}{ll}
            \ov(\alvec)\cdot\gapr\cdot\om&\mbox{if } \ga_1=\aln\\[2mm]
            \ov(\alvec^\frown\ga\cdot\om)&\mbox{otherwise.} 
            \end{array}\right.\]
If $\ga<\aln$ we have $(\kvalga+\dpval(\ga))\cdot\om<\ov(\alvec)$. 
\item For $\ga\in\domkval-(\Ez\cup\{0\})$ we have \[\dpval(\ga)<\kvalga.\] 
\item For $\ga\in\Ez^{>\aln}$ such that $\muga<\ga$ we have \[\dpval(\ga)<\ov(\alvec^\frown\ga)\cdot\muga\cdot\om.\]
\item For $\ga\in\domkval\cap\Ez^{>\aln}$ we have 
\[\kvalga\cdot\om\le\dpval(\ga) \quad\mbox{ and }\quad \dpval(\ga)\cdot\om=_\NF\ov(\homega(\alvecga))\cdot\om.\]
\item Let $\ga\in\domnuval\cap\Hz$,  $\ga=_\MNF\ga_1\cdot\ldots\cdot\ga_k$. We have
\[(\nuvalga+\ka^\alvec_{\varrho^\aln_\ga}+\dpval(\varrho^\aln_\ga))\cdot\om=\left\{\begin{array}{ll}
            \ov(\alvec)\cdot\ga\cdot\om&\mbox{if } \ga_1\le\aln\\[2mm]
            \ov(\homega(\alvec^\frown\ga))\cdot\om&\mbox{if } \ga\in\Ez^{>\aln}\\[2mm]
            \ov(\alvec^\frown\ga)\cdot\om&\mbox{otherwise.} 
            \end{array}\right.\]
\end{enumerate}
\end{lem}
{\bf Proof.} The lemma is shown by simultaneous induction on $\lalvec(\ga)$ over all parts.
\\[2mm]
{\bf Ad 1.} The claim is immediate if $\ga=\aln$, and if $\ga_1=\aln$ and $k>1$, we have $\kvalga=\ov(\alvec)\cdot\gapr$ and by part 2
the claim follows. Now assume that $\ga_1>\aln$, whence $\kvalga=\ov(\alvecga)$.
The case $\ga\not\in\Ez$ is handled again by part 2. If $\ga\in\Ez$, we apply part 4 to see that
\[(\kvalga+\dpval(\ga))\cdot\om=\dpval(\ga)\cdot\om=\ov(\homega(\alvecga))\cdot\om,\]
and since $\muga\ge\om$, the latter is equal to $\ov(\alvec^\frown\ga\cdot\om)$.

Now consider the situation where $\ga<\aln$. Let $i\in[0,\ldots,n-1]$ be maximal such that $\ali<\ga$.
The same argument as above yields the corresponding claim for $\alvec_{\restriction_i}$ instead of $\alvec$, and by the 
monotonicity of $\ov$ we see that the resulting ordinal is strictly below $\ov(\alvec_{\restriction_{i+1}})\le\ov(\alvec)$.

Note that in the case $\ga\cdot\om\in\domkval$ we have $(\kvalga+\dpval(\ga))\cdot\om=\kval_{\ga\cdot\om}$ as a direct consequence
of the definitions.\\[2mm]
{\bf Ad 2.} We may assume that $\ga>\aln$ (otherwise replace $n$ by the suitable $i<n$ and $\alvec$ by $\alvec_{\restriction_i}$).
Set $\gapr:=\log((1/\aln)\cdot\ga)=_\ANF\ga_1+\ldots+\ga_m$, so that $\ga=\aln\cdot\om^{\ga_1+\ldots+\ga_m}$ and
$\dpval(\ga)=\sum^m_{i=1}(\kval_{\ga_i}+\dpval(\ga_i))$. According to the definition, $\kvalga$ is either 
$\ov(\alvec)\cdot\om^{\ga_1+\ldots+\ga_m}$ if $\ga_1\le\aln$, or $\ov(\alvecga)$ if $\ga_1>\aln$. 
In the case $\ga_i\cdot\om<\ga$ for $i=1,\ldots,m$ an application
of part 1 of the i.h.\ to the $\ga_i$ yields the claim, thanks to the monotonicity of $\ov$. 
Otherwise we must have $\ga=\ga_1\cdot\om$ where $\ga_1\in\Ez^{>\aln}$, and
applying part 1 of the i.h.\ to $\ga_1$ yields
\[\dpval(\ga)=\kval_{\ga_1}+\dpval(\ga_1)+1<(\kval_{\ga_1}+\dpval(\ga_1))\cdot\om=\ov(\alvecga)=\kvalga.\]  
{\bf Ad 3.} Let $\la_\ga=_\ANF\la_1+\ldots+\la_r$, $\la\in\{\la_1,\ldots,\la_r\}$, and note that $\la\le\ga\cdot\muga$.
Setting $\lapr:=(1/\mf(\la))\cdot\la$, we have $\lapr\le\muga$ and applying part 1 of the i.h.\ to $\la$ we obtain
$(\kval_\la+\dpval(\la))\cdot\om\le\ov(\alvecga)\cdot\muga\cdot\om$.
Since $\nu^{\alvecga}_{\muga}=\ov(\alvec^\frown(\ga,\muga))=\ov(\alvecga)\cdot\muga$, we obtain the claim.\\[2mm]
{\bf Ad 4.} The inequality is seen by a quick inspection of the respective definitions. We have $\kvalga=\ov(\alvecga)$,
and since $\muga\ge\om$ we obtain
\[\ov(\alvecga)\cdot\om=\ov(\alvec^\frown(\ga,\om))\le\ov(\alvec^\frown(\ga,\muga))\le\dpval(\ga).\]
In order to verify the claimed equation, note that \[\homega(\alvecga)=\alvec^\frown\sk_\om(\ga),\]
where $\sk_\om(\ga)=(\de_1,\ldots,\de_{l+1})$ consists of a maximal strictly increasing chain 
$\devec:=(\de_1,\ldots,\de_l)=(\ga,\muga,\mu_{\muga},\ldots)$ of $\Ez$-numbers and $\de_{l+1}=\mu_{\de_l}\not\in\Ez^{>\de_l}$.
We have $\la_{\de_i}=\varrho_{\mu_{\de_i}}+\ze_{\de_i}$ where $\ze_{\de_i}<\de_i$ and, for $i<l$, 
$\varrho_{\mu_{\de_i}}=\mu_{\de_i}=\de_{i+1}\in\Ez^{>\de_i}$. Applying the i.h.\ to (the additive decompositions) of these
terms $\la_{\de_i}$ we obtain
\[\dpval(\ga)\cdot\om=\dpf_{\alvec^\frown\devec_{\restriction_{l-1}}}(\de_l)\cdot\om=
(\nu^{\alvec^\frown\devec}_{\mu_{\de_l}}+\ka^{\alvec^\frown\devec}_{\la_{\de_l}}+\dpf_{\alvec^\frown\devec}(\la_{\de_l}))\cdot\om,\]
and consider the additive decomposition of the term $\la_{\de_l}$. The components below $\de_l$ from $\ze_{\de_l}$ are easily handled 
using the inequality of part 1 of the i.h., while for $\mu_{\de_l}=_\MNF\mu_1\cdot\ldots\cdot\mu_j$ we have
$\varrho_{\mu_{\de_l}}\le\de_l\cdot\log(\mu_{\de_l})=\de_l\cdot(\log(\mu_1)+\ldots+\log(\mu_j))$ (where the only possible difference is
$\de_l$) and consider the summands separately.
Let $\mu\in\{\mu_1,\ldots,\mu_j\}$.\\[2mm]
{\bf Case 1:} $\varrho_{\mu_{\de_l}}\in\Ez^{>\de_l}$. Let $\log(\mu_{\de_l})=\la+k$ where $\la\in\Lim\cup\{0\}$ and $k<\om$.
We must have $\chi^{\de_l}(\la)=1$, since otherwise $\varrho_{\mu_{\de_l}}=\mu_{\de_l}\in\Ez^{>\de_l}$, which would contradict the 
maximality of the length of $\devec$.
It follows that $k=1$, hence $\mu_{\de_l}=\la\cdot\om$, $\la=\varrho_{\mu_{\de_l}}\in\Ez^{>\de_l}$, and applying the i.h.\ to $\la$
yields \[(\ka^{\alvec^\frown\devec}_{\la}+\dpf_{\alvec^\frown\devec}(\la))\cdot\om=\ov(\homega(\alvec^\frown\devec^\frown\la))\cdot\om=
\ov(\alvec^\frown\devec^\frown\la\cdot\om)=\nu^{\alvec^\frown\devec}_{\mu_{\de_l}},\]
whence $\dpval(\ga)\cdot\om=\ov(\homega(\alvec^\frown\ga))\cdot\om$ as claimed.
\\[2mm]
{\bf Case 2:} Otherwise.\\[2mm]
{\bf Subcase 2.1:} $\mu<\de_l$. Then applying part 2 of the i.h.\ to $\de_l\cdot\log(\mu)$ we see that
\[\dpf_{\alvec^\frown\devec}(\de_l\cdot\log(\mu))<\ka^{\alvec^\frown\devec}_{\de_l\cdot\log(\mu)}=\ov(\alvec^\frown\devec)\cdot\log(\mu).\] 
{\bf Subcase 2.2:} $\mu=\de_l$. We calculate $\ka^{\alvec^\frown\devec}_{\de_l^2}=\ov(\alvec^\frown\devec)\cdot\de_l$ and 
$\dpf_{\alvec^\frown\devec}(\de_l^2)=\ov(\alvec^\frown\devec)$.
\\[2mm]
{\bf Subcase 2.3:} $\mu>\de_l$. 
\\[2mm]
{\bf 2.3.1:} $\mu\not\in\Ez^{>\de_l}$. Then $\de_l<\de_l\cdot\log(\mu)\not\in\Ez^{>\de_l}$, hence by the i.h., applied to $\de_l\cdot\log(\mu)$, 
which is a summand of $\la_{\de_l}$, 
\[\dpf_{\alvec^\frown\devec}(\de_l\cdot\log(\mu))<\ka^{\alvec^\frown\devec}_{\de_l\cdot\log(\mu)}\le\nu^{\alvec^\frown\devec}_\mu\le\nu^{\alvec^\frown\devec}_{\mu_{\de_l}}.\]
{\bf 2.3.2:} $\mu\in\Ez^{>\de_l}$. Then we have $\de_l\cdot\log(\mu)=\mu$, $\mu\cdot\om\le\mu_{\de_l}$, and applying the i.h.\ to $\mu$
\[\dpf_{\alvec^\frown\devec}(\mu)\cdot\om=\ov(\homega(\alvec^\frown\devec^\frown\mu))\cdot\om=\nu^{\alvec^\frown\devec}_{\mu\cdot\om}\le\nu^{\alvec^\frown\devec}_{\mu_{\de_l}}.\]  
{\bf Ad 5.} We have $\varrho_\ga\le\aln\cdot(\log(\ga_1)+\ldots+\log(\ga_k))$. \\[2mm]
{\bf Case 1:} $\ga\in\Ez^{>\aln}$. Here we have $\varrho_\ga=\ga$ and $\nuvalga=\kvalga<\dpval(\ga)=\nu^\alvecga_{\muga}+\ka^\alvecga_{\la_\ga}+\dpf_\alvecga(\la_\ga)$,
and by part 4 we have $\dpval(\ga)\cdot\om=\ov(\homega(\alvecga))\cdot\om$.
\\[2mm]
{\bf Case 2:} $\ga_1\le\aln$. Here we argue similarly as in the proof of part 4, case 2. However, the access to the i.h.\ is different.
 Clearly, $\kval_{\aln\cdot\log(\ga_i)}=\ov(\alvec)\cdot\log(\ga_i)$. For given $i$, let $\log(\log(\ga_i))=_\ANF\xi_1+\ldots+\xi_s$. In the case $\ga_i=\aln$ we have $\dpval(\al_n^2)=\ov(\alvec)$. Now assume that $\ga_i<\aln$. An application of part 1 of the i.h.\ to $\xi_j$ yields $(\kval_{\xi_j}+\dpval(\xi_j))\cdot\om<\ov(\alvec)$ 
for $j=1,\ldots,s$.
Therefore \[(\kval_{\aln\cdot\log(\ga_i)}+\dpval(\aln\cdot\log(\ga_i)))\cdot\om=\kval_{\aln\cdot\log(\ga_i)}\cdot\om=\ov(\alvec)\cdot\log(\ga_i)\cdot\om.\]
These considerations show that we obtain \[(\nuvalga+\ka^\alvec_{\varrho^\aln_\ga}+\dpval(\varrho^\aln_\ga))\cdot\om=\nuvalga\cdot\om=\ov(\alvec)\cdot\ga\cdot\om.\]
{\bf Case 3:} Otherwise. 
\\[2mm]
{\bf Subcase 3.1:} $\varrho_\ga\in\Ez^{>\aln}$. Since $\ga\not\in\Ez^{>\aln}$ we have $\log(\ga)=\la+1$ where $\la\in\Ez^{>\aln}$ and $\chi^\aln(\la)=1$.
Hence $\ga=\la\cdot\om$ and $\varrho_\ga=\la$, for which part 4 of the i.h. yields 
\[(\kval_\la+\dpval(\la))\cdot\om=\dpval(\la)\cdot\om=\ov(\homega(\alvec^\frown\la))\cdot\om=\ov(\alvec^\frown\ga),\]
implying the claim.\\[2mm]
{\bf Subcase 3.2:} $\varrho_\ga\not\in\Ez^{>\aln}$. Here we extend the argumentation from case 2, where the situation $\ga_i\le\aln$ has been resolved.
In the case $\ga_i\in\Ez^{>\aln}$ we have $\ga_i\cdot\om\le\ga$ and apply part 4 of the i.h.\ to $\ga_i$ to obtain
\[(\kval_{\aln\cdot\log(\ga_i)}+\dpval(\aln\cdot\log(\ga_i)))\cdot\om=\dpval(\ga_i)\cdot\om=\ov(\alvec^\frown\ga_i\cdot\om)\le\ov(\alvecga).\]
We are left with the cases where $\ga_i\in\Mz^{>\aln}-\Ez$. Writing $\log(\log(\ga_i))=_\ANF\xi_1+\ldots+\xi_s$, which resides in $(\aln,\log(\ga_i))$, we
have $\dpval(\aln\cdot\log(\ga_i))=\sum^s_{j=1}(\kval_{\xi_j}+\dpval(\xi_j))$, where for each $j$ part 1 of the i.h.\ applied to $\xi_j$ together with the
monotonicity of $\kval$ on additive principal arguments yield \[(\kval_{\xi_j}+\dpval(\xi_j))\cdot\om=\kval_{\xi_j\cdot\om}\le\kval_{\aln\cdot\log(\ga_i)}<\nuvalga,\]
and we conclude as in case 2.\qed

\begin{cor}\label{kappanuhzcor} Let $\alvec\in\RS$. We have
\begin{enumerate}
\item $\kval_{\ga\cdot\om}=(\kvalga+\dpval(\ga))\cdot\om$ for $\ga\in\Hz$ such that $\ga\cdot\om\in\domkval$.
\item $\nuval_{\ga\cdot\om}=(\nuvalga+\ka^\alvec_{\varrho^\aln_\ga}+\dpval(\varrho^\aln_\ga))\cdot\om$ for $\ga\in\Hz$ such that $\ga\cdot\om\in\domnuval$.
\end{enumerate}
$\kval$ and for $\alvec\not=()$ also $\nuval$ are strictly monotonically increasing and continuous.
\end{cor}
{\bf Proof.} Parts 1 and 2 follow from Definitions \ref{odef} and \ref{kappanuprincipals} using parts 1 and
5 of Lemma \ref{kdpmainlem}, respectively. 
In order to see general monotonicity and continuity we can build upon Corollary \ref{kappanucor}.
The missing argument is as follows.
For any $\be\in\Hz$ in the respective domain and any $\ga=_\ANF\ga_1+\ldots+\ga_m<\be$
we have $\kvalga<\kvalbe$ using part 1, and $\nuvalga<\nuvalbe$ using part 2, since $\ga_i\cdot\om\le\be$
for $i=1,\ldots,m$.\qed

\begin{theo}\label{agreementthm}
Let $\alvec=(\ale,\ldots,\aln)\in\RS$, $n\ge 0$, and set $\al_0:=1$.
For $\be\in\Hz$ let $\de:=(1/\bebar)\cdot\be$, so that $\be=_\NF\bebar\cdot\de$ if $\be\not\in\Mz$.
\begin{enumerate}
\item For all $\be\in\domkval\cap\Hz^{>\aln}$ we have \[\kvalbe=\kval_{\bebar+1}\cdot\de.\]
\item For all $\be\in\domnuval\cap\Hz^{>\aln}$ (where $n>0$) we have \[\nuvalbe=\nuval_{\bebar+1}\cdot\de.\]
\end{enumerate}
Hence, the definitions of $\kappa,\nu$, and $\dpf$ given in \cite{CWc} and here fully agree.
\end{theo}
{\bf Proof.} We rely on the monotonicity of $\ov$. Note that Corollary \ref{kappanuhzcor} has already shown the theorem for $\be$ of the form $\ga\cdot\om$,
i.e.\ successors of additive principal numbers. Let $\be=_\MNF\be_1\cdot\ldots\cdot\be_k\in\Hz^{>\aln}$.\\[2mm]
{\bf Case 1:} $k=1$. Then $\de=\be$, and setting $n_0:=n_0(\alvecbe)$ and $\ga:=\ga(\alvecbe)$ according to Definition \ref{odef}, we have 
\[\kvalbe=\nuvalbe=\ov(\alvecbe)=\left\{\begin{array}{ll}
            \be&\mbox{if } n_0=0\\[2mm]
            \ov(\alvecpr)\cdot\be&\mbox{if }n_0>0,
            \end{array}\right.\]
where $\alvecpr:=\hbe({\alvec_{\restriction_{n_0-1}}}^\frown\ga)$, and parts 1 and 5 of Lemma \ref{kdpmainlem} yield
\[\kval_{\bebar+1}\cdot\be=\nuval_{\bebar+1}\cdot\be=\left\{\begin{array}{ll}
            \ov(\alvec)\cdot\be&\mbox{if } \bebar=\aln\\[2mm]
            \ov(\homega(\alvec^\frown\bebar))\cdot\be&\mbox{if }\bebar>\aln.
            \end{array}\right.\]     
If $n_0=0$ the claim is immediate since $\ov(\alvec)\cdot\be=\be$ if $\bebar=\aln$ and 
$1<\ov(\homega(\alvec^\frown\bebar))<\ov(\alvecbe)=\be$ if $\bebar>\aln$. Now assume that $n_0>0$.\\[2mm]
{\bf Subcase 1.1:} $\bebar=\aln$. This implies $n_0\le n$ and therefore $\alvecpr \klex\alvec\klex\alvecbe$.
By the monotonicity of $\ov$ we obtain
\[\ov(\alvecpr )<\ov(\alvecbe)=\ov(\alvecpr )\cdot\be,\]
which implies the claim since $\be\in\Mz$.\\[2mm]
{\bf Subcase 1.2:} $\bebar>\aln$. This implies $\bebar\in\Ez^{>\aln}$ and $\kval_{\bebar+1}\cdot\be=\nuval_{\bebar+1}\cdot\be=\ov(\homega(\alvec^\frown\bebar))\cdot\be$,
where $\homega(\alvec^\frown\bebar)\klex\alvecbe$.\\[2mm]
{\bf 1.2.1:} $n_0\le n$. Then we obtain
$\alvecpr \klex\alvec\klex\homega(\alvec^\frown\bebar)\klex\alvecbe$ and hence
\[\ov(\alvecpr )<\ov(\homega(\alvec^\frown\bebar))<\ov(\alvecbe)=\ov(\alvecpr )\cdot\be,\]
which implies the claim.\\[2mm]
{\bf 1.2.2:} $n_0=n+1$. Then we have $\ga\le\bebar$, $\alvecpr=\hbe(\alvecga)$, and using Corollary \ref{hgacor} it follows that 
\[\kvalbe=\nuvalbe=\ov(\alvecbe)=\ov(\alvecpr)\cdot\be>\ov(\homega(\alvec^\frown\bebar))\ge\ov(\hbe(\alvecga)),\]
which again implies the claim.\\[2mm]
{\bf Case 2:} $k>1$. Then we have $\bebar=_\MNF\be_1\cdot\ldots\cdot\be_{k-1}\ge\aln$, $\de=\be_k$, and set $\bepr:=(1/\be_1)\cdot\be$ and $\bebarpr:=(1/\be_1)\cdot\bebar$.\\[2mm]
{\bf Subcase 2.1:} $\be_1=\aln$. Using Lemma \ref{kdpmainlem} we obtain
\[\kval_{\bebar+1}\cdot\de=\ov(\alvec)\cdot\bepr=\kvalbe\]
and
\[\nuval_{\bebar+1}\cdot\de=\ov(\alvec)\cdot\be=\nuvalbe.\]
{\bf Subcase 2.2:} $\be_1>\aln$. Then we have $\kvalbe=\nuvalbe$.
\\[2mm]
{\bf 2.2.1:} $\bebar\not\in\Ez^{>\aln}$. Then by the involved definitions 
\[\kval_{\bebar+1}\cdot\de=\nuval_{\bebar+1}\cdot\de=\ov(\alvec^\frown\bebar)\cdot\de=\nuvalbe=\kvalbe.\]
{\bf 2.2.2:} $\bebar\in\Ez^{>\aln}$. This implies $k=2$, $\de=\be_2$, and we see that
\[\kval_{\bebar+1}\cdot\de=\nuval_{\bebar+1}\cdot\de=\ov(\homega(\alvec^\frown\bebar))\cdot\de.\]
In the case $\de>\mu_{\be_1}$ we have $\hop_{\de}(\alvec^\frown\bebar)=\alvec^\frown\bebar$ and hence obtain uniformly
\[\kvalbe=\nuvalbe=\ov(\hop_{\de}(\alvec^\frown\bebar))\cdot\de.\]
By Corollary \ref{hgalem} we have $\hop_{\de}(\alvec^\frown\bebar)\kglex\homega(\alvec^\frown\bebar)$, hence
\[\ov(\hop_{\de}(\alvec^\frown\bebar))\le\ov(\homega(\alvec^\frown\bebar))<\ov(\alvecbe)=\ov(\hop_{\de}(\alvec^\frown\bebar))\cdot\de,\]
which implies the claim since $\de\in\Mz$.
\qed

\begin{lem}\label{hatmainlem}
Let $\alvec=(\ale,\ldots,\aln)\in\RS$, $n>0$. 
\begin{enumerate}
\item For all $\be$ such that $\alvecbe\in\TS$ we have 
\[\ov(\alvecbe)<\ov(\alvec)\cdot\widehat{\aln}.\]
\item For all $\ga$ such that $\alvecga\in\RS$ we have 
\[\ov(\homega(\alvecga))<\ov(\alvecga)\cdot\gahat.\]
\end{enumerate}
\end{lem}
{\bf Proof.} We prove the lemma by simultaneous induction on $\lSeq(\alvecbe)$ and $\lSeq(\homega(\alvecga))$, respectively.\\[2mm]
{\bf Ad 1.}
Let $\be=_\MNF\be_1\cdot\ldots\cdot\be_k$ and $\bepr:=(1/\be_1)\cdot\be$.
Note that $\be\le\mualn<\alhat$.
\\[2mm]
{\bf Case 1:} $\be_1\le\aln$. Immediate, since $\ov(\alvecbe)=\ov(\alvec)\cdot\be$.
\\[2mm]
{\bf Case 2:} $k>1$ where $\be_1\in\Ez^{>\aln}$ and $\be_2\le\mu_{\be_1}$. Then we have
$\ov(\alvecbe)=\ov(\hop_{\be_2}(\alvec^\frown\be_1))\cdot\bepr$ and apply the i.h.\ to $\alvec^\frown\be_1$,
which clearly satisfies $\lSeq(\alvec^\frown\be_1)\klex\lSeq(\alvecbe)$.
By Corollaries \ref{hgacor} and \ref{ocontcor} we have 
$\ov(\hop_{\be_2}(\alvec^\frown\be_1))\le\ov(\homega(\alvec^\frown\be_1))$,
and by the i.h., parts 1 (for $\alvec^\frown\be_1$) and 2 (for $\homega(\alvec^\frown\be_1)$) we obtain
\[\ov(\homega(\alvec^\frown\be_1))<\ov(\alvec^\frown\be_1)\cdot\widehat{\be_1}<\ov(\alvec)\cdot\widehat{\aln},\]
where we have used that $\widehat{\be_1}\le\widehat{\aln}$ according to Lemma 3.17 of \cite{CWc}.
This implies the desired inequality.
\\[2mm]
{\bf Case 3:} Otherwise. Let $n_0:=n_0(\alvecbe)$ and $\ga:=\ga(\alvecbe)$ according to Definition \ref{odef}.
\\[2mm]
{\bf Subcase 3.1:} $n_0=0$. Immediate.
\\[2mm]
{\bf Subcase 3.2:} $n_0>0$. By definition we have 
$\ov(\alvecbe)=\ov(\hop_{\be_1}({\alvec_{\restriction_{n_0-1}}}^\frown\ga))\cdot\be$.
\\[2mm]
{\bf 3.2.1:} $n_0\le n$. The monotonicity of $\ov$ then yields 
$\ov(\alvecbe)\le\ov(\alvec)\cdot\be<\ov(\alvec)\cdot\widehat{\aln}$.
\\[2mm]
{\bf 3.2.2:} $n_0=n+1$. Then $\ga$ is the immediate predecessor of $\be$ in $\trs^\aln(\be)$.
We apply the i.h.\ for parts 1 (to $\alvecga$) and 2 (to $\homega(\alvecga)$) and argue as in Case 2 to see that
$\ov(\alvecbe)=\ov(\hop_{\be_1}(\alvecga))\cdot\be<\ov(\alvec)\cdot\widehat{\aln}$.\\[2mm]
{\bf Ad 2.} We have $\homega(\alvecga)=\alvec^\frown\sk_\om(\ga)$, and setting
$\sk_\om(\ga)=:(\de_1,\ldots,\de_{l+1})$ we obtain a strictly increasing sequence of $\Ez$-numbers
$\devec:=(\de_1,\ldots,\de_l)=(\ga,\muga,\mu_{\muga},\ldots)$ such that 
$\de_{l+1}=\mu_{\de_l}\not\in\Ez^{>\de_l}$. For $i=1,\ldots,l$ we have
\[\homega(\alvec^\frown\devec_{\restriction_i})=\alvec^\frown\devec^\frown\de_{l+1},\]
\[\lSeq(\alvec^\frown\devec_{\restriction_i})\kglex\lSeq(\homega(\alvecga)),\]
and the i.h., part 1 (up to $\homega(\alvecga)$), yields
\[\ov(\alvec^\frown\devec_{\restriction_{i+1}})<\ov(\alvec^\frown\devec_{\restriction_i})\cdot\widehat{\de_i}.\]
Appealing to Lemma 3.17 of \cite{CWc} we obtain
\[\widehat{\de_l}\le\widehat{\de_{l-1}}\le\ldots\le\widehat{\de_1}=\gahat,\]
and finally conclude that $\ov(\homega(\alvecga))<\ov(\alvecga)\cdot\gahat$.
\qed

\begin{cor}\label{kdpnuestimcor}
For all $\alvecga\in\RS$ the ordinal $\ov(\alvecga)\cdot\gahat$ is a strict upper bound of  
\[\Image(\kappa^{\alvecga}),\: \Image(\nu^{\alvecga}),\: \dpval(\ga),\: \mbox{and }
\nu^{\alvecga}_{\mu_\ga}+\kappa^{\alvecga}_{\la_\ga}+\dpf_{\alvecga}(\la_\ga).\]
\end{cor}
{\bf Proof.} This directly follows from Lemmas \ref{kdpmainlem} and \ref{hatmainlem}.
\qed

\section{Revisiting tracking chains}\label{revisitsec}

\subsection{Preliminary remarks}\label{prelimsubsec}
Our preparations in the previous sections are almost sufficient to demonstrate that the characterization of 
$\Ctwo$ provided in Section 7 of \cite{CWc} is elementary recursive. We first provide a brief
argumentation based on the previous sections showing that the structure $\Ctwo$ is elementary recursive.
In the following subsections we will elaborate on the characterization of $\leo$ and $\letwo$
within $\Ctwo$, further illuminating the structure.

In Section 5 of \cite{CWc} the termination of the process of maximal extension (see Definition 5.2 of \cite{CWc})
is seen when applying the $\ltvec$-measure from the second step on, as clause 2.3.1 of Def.\ 5.2 in \cite{CWc} can
only be applied at the beginning of the process of maximal extension.
Lemma 5.4, part a), of \cite{CWc} is not needed in full generality, the $\ltvec$-measure suffices, cf.\ 
Lemma 5.5 of \cite{CWc}.
The proof of Lemma 5.12 of \cite{CWc}, parts a), b), and c), actually proceeds by induction on the number
of 1-step extensions, an ``induction on $\cspr(\alvec)$ along $\klex$'' is not needed.

In Section 6 of \cite{CWc} the proof of Lemma 6.2 actually proceeds by induction on the length of the additive decomposition of $\al$. Definition 6.1 of \cite{CWc}, which assigns to each $\al<\oneinf$ its unique tracking chain $\tc(\al)$, involves the evaluation function $\ov$ 
(in the guise of $\tilde{\cdot}$, cf.\ Definition 5.9 of \cite{CWc}), which we have shown to be elementary recursive.

\subsection{Pre-closed and spanning sets of tracking chains}
For the formal definition of tracking chains recall Definition 5.1 of \cite{CWc}. We will rely on its detailed terminology,
including the notation $\alvec\sub\bevec$ when $\alvec$ is an initial chain of $\bevec$.

The following definitions of pre-closed and spanning sets of tracking chains provide a generalization of the notion of maximal 
extension, denoted by $\me$, cf.\ Definition 5.2 of \cite{CWc}. 

\begin{defi}[Pre-closedness]\label{precldefi} Let $M\finsub\TC$. $M$ is \emph{pre-closed} if and only if $M$ 
\begin{enumerate}
\item is \emph{closed under initial chains:} if $\alvec\in M$ and $(i,j)\in\dom(\alvec)$ then $\alvec_{\restriction_{(i,j)}}\in M$,
\item is \emph{$\nu$-index closed:} if $\alvec\in M$, $m_n>1$, $\alcp{n,m_n}=_\ANF\xi_1+\ldots+\xi_k$ then
$\alvec[\mu_{\taupr}], \alvec[\xi_1+\ldots+\xi_l]\in M$ for $1\le l\le k$,
\item \emph{unfolds minor $\letwo$-components:} if $\alvec\in M$, $m_n>1$, and $\tau<\mu_\taupr$ then: 
\begin{enumerate}
\item[3.1.] ${\alvec_{\restriction_{n-1}}}^\frown(\alcp{n,1},\ldots,\alcp{n,m_n},\mu_\tau)\in M$ in the case $\tau\in\Ez^{>\taupr}$, and
\item[3.2.] otherwise $\alvec^\frown(\varrho^\taupr_\tau)\in M$, provided that $\varrho^\taupr_\tau>0$, 
\end{enumerate}
\item is \emph{$\ka$-index closed:} if $\alvec\in M$, $m_n=1$, and $\alcp{n,1}=_\ANF\xi_1+\ldots+\xi_k$, then:
\begin{enumerate}
\item[4.1.] if $m_{n-1}>1$ and $\xi_1=\taucp{n-1,m_{n-1}}\in\Ez^{>\taucp{n-1,m_{n-1}-1}}$ then 
${\alvec_{\restriction_{n-2}}}^\frown(\alcp{n-1,1},\ldots,\alcp{n-1,m_{n-1}},\mu_{\xi_1})\in M$,
else ${\alvec_{\restriction_{n-1}}}^\frown(\xi_1)\in M$, and
\item[4.2.] ${\alvec_{\restriction_{n-1}}}^\frown(\xi_1+\ldots+\xi_l)\in M$ for $l=2,\ldots,k$,
\end{enumerate}
\item \emph{maximizes $\me$-$\mu$-chains:} if $\alvec\in M$, $m_n\ge 1$, and $\tau\in\Ez^{>\taupr}$, then:
\begin{enumerate}
\item[5.1.] if $m_n=1$ then ${\alvec_{\restriction_{n-1}}}^\frown(\alcp{n,1},\mu_\tau)\in M$, and
\item[5.2.] if $m_n>1$ and $\tau=\mu_\taupr=\la_\taupr$ then ${\alvec_{\restriction_{n-1}}}^\frown(\alcp{n,1}\ldots,\alcp{n,m_n},\mu_\tau)\in M$. 
\end{enumerate}
\end{enumerate}
\end{defi}

\noindent{\bf Remark.}
Pre-closure of some $M\finsub\TC$ is obtained by closing under clauses 1 -- 5 in this order once, hence finite: 
in clause 5 note 
that $\mu$-chains are finite since the $\htarg{}$-measure of terms strictly decreases with each application of $\mu$. Note further that intermediate indices
are of the form $\la_\taupr$, whence we have a decreasing $\operatorname{l}$-measure according to inequality \ref{lalen} in the remark following Definition \ref{Ttauvec}.

\begin{defi}[Spanning sets of tracking chains]\label{spanningdefi} 
$M\finsub\TC$ is \emph{spanning} if and only if it is pre-closed and closed under
\begin{enumerate}
\item[6.] \emph{unfolding of $\leo$-components:} for $\alvec\in M$, if $m_n=1$ and $\tau\not\in\Ez^{\ge\taupr}$ 
(i.e.\ $\tau=\taucp{n,m_n}\not\in\Ezone$, $\taupr=\taunstar$), let
\[\log((1/\taupr)\cdot\tau)=_\ANF\xi_1+\ldots+\xi_k,\]
if otherwise $m_n>1$ and $\tau=\mu_\taupr$ such that $\tau<\la_\taupr$ in the case $\tau\in\Ez^{>\taupr}$, let
\[\la_\taupr=_\ANF\xi_1+\ldots+\xi_k.\]
Set $\xi:=\xi_1+\ldots+\xi_k$, unless $\xi>0$ and $\alvec^\frown(\xi_1+\ldots+\xi_k)\not\in\TC$,
\footnote{This is the case if clause 6 of Def.\ 5.1 of \cite{CWc} does not hold.} 
in which case we set $\xi:=\xi_1+\ldots+\xi_{k-1}$.
Suppose that $\xi>0$. Let $\alvecpl$ denote the vector $\{\alvec^\frown(\xi)\}$ if this is a tracking chain, 
or otherwise the vector ${\alvec_{\restriction_{n-1}}}^\frown(\alcp{n,1},\ldots,\alcp{n,m_n},\mu_\taucp{n,m_n})$.
\footnote{This case distinction, due to clause 5 of Def.\ 5.1 of \cite{CWc}, is missing in \cite{W17}.}
Then the closure of $\{\alvecpl\}$ under clauses 4 and 5 is contained in $M$.
\end{enumerate}
\end{defi}

\noindent{\bf Remark.}
Closure of some $M\finsub\TC$ under clauses 1 -- 6 is a finite process since pre-closure is finite and since 
the $\ka$-indices added in clause 6 strictly decrease in $\operatorname{l}$-measure.


Semantically, the above notion of spanning sets of tracking chains and closure under clauses 1 -- 6 leaves 
some redundancy in the form that certain
$\ka$-indices could be omitted. This will be adressed elsewhere, since the current formulation is advantageous
for technical reasons.

\begin{defi}[Relativization]\label{reldefi} Let $\alvec\in\TC\cup\{()\}$ and $M\finsub\TC$ be a set of tracking chains
that properly extend $\alvec$.
$M$ is \emph{pre-closed above $\alvec$} if and only if it is pre-closed 
with the modification that clauses 1 -- 5 only apply when the respective resulting tracking chains $\bevec$ 
properly extend $\alvec$.
$M$ is \emph{weakly spanning above $\alvec$} if and only if $M$ is pre-closed above $\alvec$ and closed under
clause 6. 
\end{defi}

\begin{lem}\label{meclosurelem} If $M$ is spanning (weakly spanning above some $\alvec$), then it is closed under $\me$
(closed under $\me$ for proper extensions of $\alvec$).
\end{lem}
{\bf Proof.} This follows directly from the definitions involved. \qed

\subsection{Characterizing $\leo$ and $\letwo$ in $\Ctwo$}
The purpose of this section is to provide a detailed picture of the restriction of $\Rtwo$ to $\oneinf$ on the basis of the 
results of \cite{CWc} and to conclude with the extraction of an elementary recursive arithmetical characterization of this 
structure given in terms of tracking chains which we will refer to as $\Ctwo$. 
 
We begin with a few observations that follow from the results in Section 7 of \cite{CWc} and explain the concept of tracking chains.
The evaluations of all initial chains of some tracking chain $\alvec\in\TC$ form a $\lo$-chain.
Evaluations of initial chains $\alvec_{\restriction_{i,j}}$ where $(i,j)\in\dom(\alvec)$ and $j=2,\ldots,m_i$ with fixed index
$i$ form $\ktwo$-chains. 
Recall that indices $\alcp{i,j}$ are $\ka$-indices for $j=1$ and $\nu$-indices otherwise, cf.\ Definitions 5.1 and 5.9 of \cite{CWc}.

According to Theorem 7.9 of \cite{CWc}, an ordinal $\al<\oneinf$ is $\leo$-minimal if and only if its tracking chain consists of 
a single $\ka$-index, i.e.\ if its tracking chain $\alvec$ satisfies $(n,m_n)=(1,1)$. Clearly, the least $\leo$-predecessor of
any ordinal $\al<\oneinf$ with tracking chain $\alvec$ is $\ov(\alvec_{\restriction_{1,1}})=\ka_\alcp{1,1}$.
According to Corollary 7.11 of \cite{CWc} the ordinal $\oneinf$ is $\leo$-minimal.
An ordinal $\al>0$ has a non-trivial $\leo$-reach if and only if $\taucp{n,1}>\taunstar$, hence in particular when $m_n>1$, 
cf.\ condition 2 in Definition 5.1 of \cite{CWc}. 

We now turn to a characterization of the \emph{greatest immediate $\leo$-successor}, $\gs(\al)$, 
of an ordinal $\al<\oneinf$ with tracking chain $\alvec$.
Recall the notations $\rho_i$ and $\alvec[\xi]$ from Definition 5.1 of \cite{CWc}. The largest $\al$-$\leo$-minimal ordinal is the root 
of the $\la$th $\al$-$\leo$-component for $\la:=\rho_n\minusp1$. Therefore, if $\al$ has a non-trivial $\leo$-reach, its greatest immediate 
$\leo$-successor $\gs(\al)$ has the tracking chain $\alvec^\frown(\la)$, \emph{unless} we either
have $\taucp{n,m_n}<\mu_{\taunpr}\andsp\chi^\taunpr(\taucp{n,m_n})=0$, where 
$\tc(\gs(\al))=\alvec[\alcp{n,m_n}+1]$, or $\alvec^\frown(\la)$ is in conflict with either condition 5 of Definition 5.1 of \cite{CWc},
in which case we have $\tc(\gs(\al))={\alvec_{\restriction_{n-1}}}^\frown(\alcp{n,1},\ldots,\alcp{n,m_n},1)$, or condition 6
of Definition 5.1 of \cite{CWc}, in which case we have $\tc(\gs(\al))=\alvec_{\restriction_{i,j+1}}[\alcp{i,j+1}+1]$.
\footnote{This condition is missing in \cite{W17}.}
In case $\al$ does not have any $\lo$-successor, we set $\gs(\al):=\al$.

$\al$ is $\letwo$-minimal if and only if for its tracking chain $\alvec$ we have $m_n\le2$ and $\taunstar=1$,
and $\al$ has a non-trivial $\letwo$-reach if and only if $m_n>1$ and $\taucp{n,m_n}>1$. Note that any $\al\in\On$ with
a non-trivial $\letwo$-reach is the proper supremum of its $\lo$-predecessors, hence $\oneinf$ does not possess any $\ktwo$-successor. 
Iterated closure under the relativized notation system $\Tt$ for $\tau=\oneinf, (\oneinf)^\infty, \ldots$ results in the infinite $\ktwo$-chain
through $\On$. Its $\lo$-root is $\oneinf$, the root of the ``master main line'' of $\Rtwo$, outside the core of $\Rtwo$, i.e.\ $\oneinf$.

According to part (a) of Theorem 7.9 of \cite{CWc} $\al$ has a greatest $\lo$-predecessor if and only if it is not
$\leo$-minimal and has a trivial $\letwo$-reach (i.e.\ does not have any $\ktwo$-successor). This is the case if 
and only if either $m_n=1$ and $n>1$, where we have $\predec_1(\al)=\ordcp{n-1,m_{n-1}}(\alvec)$, or $m_n>1$ and
$\taucp{n,m_n}=1$. In this latter case $\alcp{n,m_n}$ is of a form $\xi+1$ for some $\xi\ge 0$, and using again the notation
from Definition 5.1 of \cite{CWc} we have
$\predec_1(\al)=\ov(\alvec[\xi])$ if $\chi^\taucp{n,m_n-1}(\xi)=0$, whereas $\predec_1(\al)=\ov(\me(\alvec[\xi]))$ in
the case $\chi^\taucp{n,m_n-1}(\xi)=1$.
  
Recall Definition 7.12 of \cite{CWc}, defining for $\alvec\in\TC$ the notation $\alvec^\star$ and the index pair $\gbo(\alvec)=:(n_0,m_0)$, 
which according to Corollary 7.13 of \cite{CWc} enables us to express the $\leo$-reach $\lh(\al)$ of $\al:=\ov(\alvec)$, 
cf.\ Definition 7.7 of \cite{CWc}, by
\begin{equation}\label{lhequation}\lh(\al)=\ov(\me(\bevec^\star)),
\end{equation} 
where $\bevec:=\alvec_{\restriction_{n_0,m_0}}$, which in the case $m_0=1$ is equal to 
$\ov(\me(\bevec))=\ordcp{n_0,1}(\alvec)+\dpf_{\tauticp{n_0,0}}(\taucp{n_0,1})$
and in the case $m_0>1$ equal to $\ov(\me(\bevec[\mu_\taucp{n_0,m_0-1}]))$.
Note that if $\cml(\alvec^\star)$ does not exist we have
\[\lh(\al)=\ov(\me(\alvec^\star)),\]
and the tracking chain $\bevec$ of any ordinal $\be$ such that $\ov(\alvec^\star)\leo\be$ is then an extension of $\alvec$,
$\alvec\sub\bevec$, as will follow from Lemma \ref{subresplem}.

The relation $\leo$ can be characterized by 
\begin{equation}\label{leoequation}\al\leo\be 
\quad\aeq\quad\al\le\be\le\lh(\al),
\end{equation} 
showing that $\leo$ is a forest contained in $\le$, which \emph{respects} the ordering $\le$, i.e.\ 
if $\al\le\be\le\ga$ and $\al\leo\ga$ then $\al\leo\be$. On the basis of Lemma \ref{meclosurelem}, Equation \ref{lhequation} has the following
\begin{cor}\label{leoclosurecor} Let $M\finsub\TC$ be spanning (weakly spanning above some $\alvec\in\TC$) and $\bevec\in M$, 
$\be:=\ov(\bevec)$.
Then \[\tc(\lh(\be))\in M,\] provided that $\ov(\bevec_{\restriction_{\gbo(\bevec)}})$ is a proper extension of $\alvec$
in the case that $M$ is weakly spanning above $\alvec$.
\qed
\end{cor}
   
We now recall how to retrieve the greatest $\ktwo$-predecessor of an ordinal below $\oneinf$, if it exists, and the iteration of this procedure to obtain the maximum
chain of $\ktwo$-predecessors. Recall Definition 5.3 and Lemma 5.10 of \cite{CWc}.
Using the following proposition we can prove two other useful characterizations of the relationship $\al\letwo\be$. 
 
\begin{prop}\label{letwopredprop} Let $\al<\oneinf$ with $\tc(\al)=:\alvec$.
We define a sequence $\sivec\in\RS$ as follows.
\begin{enumerate}
\item If $m_n\le2$ and $\taunstar=1$, whence $\al$ is $\letwo$-minimal according to Theorem 7.9 of \cite{CWc}, set $\sivec:=()$. Otherwise,
\item if $m_n>2$, whence $\predec_2(\al)=\ordcp{n,m_n-1}(\alvec)$ with base $\taucp{n,m_n-2}$ according to Theorem 7.9 of \cite{CWc}, we set
$\sivec:=\cs(\alvec_{\restriction_{n,m_n-2}})$,
\item and if $m_n\le2$ and $\taunstar>1$, whence $\predec_2(\al)=\ordcp{i,j+1}(\alvec)$ with base $\taucp{i,j}$ where $(i,j):=n^\star$, again according to Theorem 7.9 of \cite{CWc}, we set
$\sivec:=\cs(\alvec_{\restriction_{i,j}})$.
\end{enumerate}
Each $\si_i$ is then of a form $\taucp{k,l}$ where $1\le l< m_k$, $1\le k \le n$. The corresponding $\ktwo$-predecessor of $\al$ is $\ordcp{k,l+1}(\alvec)=:\be_i$.
We obtain sequences $\sivec=(\si_1,\ldots,\si_r)$ and $\bevec=(\be_1,\ldots,\be_r)$ with $\be_1\ktwo\ldots\ktwo\be_r\ktwo\al$, where $r=0$ if $\al$ is $\letwo$-minimal,
so that $\predecs_2(\al)=\{\be_1,\ldots,\be_r\}$ and hence
$\be\ktwo\al$ if and only if $\be\in\predecs_2(\al)$, displaying that $\letwo$ is a forest contained in $\leo$.\qed
\end{prop}

\begin{lem}\label{letwosuclem} Let $\al,\be<\oneinf$ with tracking chains $\tc(\al)=\alvec=(\alvec_1,\ldots,\alvec_n)$, 
$\alvec_i=(\alcp{i,1},\ldots,\alcp{i,m_i})$, $1\le i\le n$, and $\tc(\be)=\bevec=(\bevec_1,\ldots,\bevec_l)$, 
$\bevec_i=(\becp{i,1},\ldots,\becp{i,k_i})$, $1\le i\le l$. 
Assume further that $\alvec\subseteq\bevec$ with
associated chain $\tauvec$ and that $m_n>1$. Set $\tau:=\taucp{n,m_n-1}$. The following are equivalent:
\begin{enumerate}
\item $\al\letwo\be$
\item $\tau\le\taucp{j,1}$ for $j=n+1,\ldots,l$
\item $\tauti\mid\be$.
\end{enumerate}
\end{lem}
{\bf Proof.} Note that $\alvec\subseteq\bevec$ is a necessary condition for $\al\letwo\be$, cf.\ Proposition \ref{letwopredprop}.
\\[2mm]
{\bf 1 $\imp$ 2:} Iterating the operation $(\cdot)^\prime$ (which for $\taucp{j,1}$ is $\taucppr{j,1}=\taujstar$) then reaches $\tau$. 
This implies $\taucp{j,1}\ge\tau$ for $j=n+1,\ldots,l$.
\\[2mm]
{\bf 2 $\imp$ 1:} Iterating the procedure to find the greatest $\letwo$-predecessor, cf.\ Proposition \ref{letwopredprop}, 
from $\be$ downward satisfies $\taujstar\ge\tau$ at each $j\in(n,l]$, where we therefore have $(n,m_n-1)\kglex j^\star$.
\\[2mm]
{\bf 2 $\imp$ 3:} Note that according to Lemma 5.10(c) of \cite{CWc} we have \[\trs(\tauti)=\cs(\alvec_{\restriction_{n,m_n-1}}).\] 
We argue by induction along $\klex$ on $(l,k_l)$. The case $k_l>1$ is trivial since $\be$ is then a multiple of $\ov(\bevec_{\restriction_{l,k_l-1}})$.
Assume now that $k_l=1$, so that $\taucp{l,1}\ge\tau$. Let $(u,v):=l^\star$, so $(n,m_n-1)\kglex(u,v)$ and $\tau_l^\star\ge\tau$. By Lemma 5.10(b) of \cite{CWc}
we have \[\tauticp{l,1}=\ka^{\tauti^\star_l}_\taucp{l,1}.\] If $(n,m_n-1)=(u,v)$ we obtain $\tauti\mid\ka^{\tauti}_\taucp{l,1}$ since $\taucp{l,1}\ge\tau$.
By the i.h.\ we have $\tauti\mid\ov(\bevec_{\restriction_{u,v}})$, so $\tauti\mid\tauticp{u,v}$, and since $\taucp{l,1}\ge\tau^\star_l$ we also have
$\tauti\mid\tauticp{l,1}$, which implies that $\tauti\mid\be$, cf.\ Definition \ref{kappanuprincipals}.
\\[2mm]
{\bf 3 $\imp$ 2:} Assume there exists a maximal $j\in\{n+1,\ldots,l\}$ such that $\taucp{j,1}<\tau$. Then we have $\taujstar\le\taucp{j,1}<\tau$, so 
$(u,v):=j^\star\klex(n,m_n-1)$. Let $\trs(\tauti)=:(\si_1,\ldots,\si_s)$ and recall Theorem \ref{thma}.
\\[2mm]
{\bf Case 1:} $\taucp{j,1}\not\in\Ez^{>\taujstar}$. Then it follows that $(j,1)=(l,k_l)$ since $j$ was chosen maximally.
In the case $\taujstar=1$ we have $\taucp{j,1}<\si_1$ and hence $\tauticp{j,1}=\ka_\taucp{j,1}<\tauti$. Otherwise, i.e.\ $\taujstar=\taucp{u,v}>1$, we
obtain $\trs(\tauticp{u,v})=\cs(\alvec_{\restriction_{u,v}})\klex\trs(\tauti)$, so \[\tauticp{j,1}=\ka^{\trs(\tauticp{u,v})}_\taucp{j,1}<\tauti.\] 
{\bf Case 2:} $\taucp{j,1}\in\Ez^{>\taujstar}$. Then $\trs(\tauticp{j,1})=\cs(\bevec_{\restriction_{j,1}})=\cs(\bevec_{\restriction_{u,v}})^\frown\taucp{j,1}\klex\trs(\tauti)$,
hence $\tauticp{j,1}<\tauti$. In the case $\taucp{l,k_l}\in\Ez^{>\taucppr{l,k_l}}$ we see that $\trs(\tauticp{j,1})\subseteq\trs(\tauticp{l,k_l})$, hence 
$\tauticp{l,k_l}<\tauti$. The other cases are easier, cf.\ Case 1.
\qed

Applying the mappings $\tc$, cf.\ Section 6 of \cite{CWc}, and $\ov$, which we verified in the present article to be elementary recursive,
cf.\ Subsection \ref{prelimsubsec},
we are now able to formulate the arithmetical characterization of $\Ctwo$.

\begin{cor}\label{elemreccharcor} The structure $\Ctwo$ is characterized elementary recursively by
\begin{enumerate}
\item $(\oneinf,\le)$ is the standard ordering of the classical notation system $\oneinf=\Targ{1}\cap\Om_1$, cf.\ \cite{W07a},
\item $\al\leo\be$ if and only if $\al\le\be\le\lh(\al)$ where $\lh$ is given by equation \ref{lhequation}, and
\item $\al\letwo\be$ if and only if $\tc(\al)\subseteq\tc(\be)$ and condition 2 of Lemma \ref{letwosuclem} holds.\qed
\end{enumerate}
\end{cor}

\begin{cor}\label{lhtwoclscor}
Let $M\finsub\TC$ be spanning (weakly spanning above some $\alvec\in\TC$). Then $M$ is closed under $\lh_2$.
\end{cor}
{\bf Proof.} This follows from Lemma \ref{meclosurelem} using Lemma \ref{letwosuclem}, cf.\ Corollaries 5.6 and 7.13 of \cite{CWc}. 
\qed

Recall Definition 5.13 of \cite{CWc} which characterizes the standard linear ordering $\le$ on $\oneinf$
by an ordering $\letc$ on the corresponding tracking chains.
We can formulate a characterization of the relation $\leo$ (below $\oneinf$) in terms of the corresponding 
tracking chains as well. This follows from an inspection of the ordering $\letc$ in combination with the above 
statements. Let $\al,\be<\oneinf$ with tracking chains $\tc(\al)=\alvec=(\alvec_1,\ldots,\alvec_n)$, 
$\alvec_i=(\alcp{i,1},\ldots,\alcp{i,m_i})$, $1\le i\le n$, and $\tc(\be)=\bevec=(\bevec_1,\ldots,\bevec_l)$, 
$\bevec_i=(\becp{i,1},\ldots,\becp{i,k_i})$, $1\le i\le l$. 
We have $\al\leo\be$ if and only if either $\alvec\subseteq\bevec$ or there exists 
$(i,j)\in\dom(\alvec)\cap\dom(\bevec)$, $j<\min\{m_i,k_i\}$, such that 
$\alvec_{\restriction_{i,j}}=\bevec_{\restriction_{i,j}}$ and $\alcp{i,j+1}<\becp{i,j+1}$,
and we either have $\chi^\taucp{i,j}(\alcp{i,j+1})=0\andsp(i,j+1)=(n,m_n)$ or 
$\chi^\taucp{i,j}(\alcp{i,j+1})=1\andsp\al\leo\ov(\me(\alvec_{\restriction_{i,j+1}}))\lo\be$.
Iterating this argument and recalling Lemma 5.5 of \cite{CWc} we obtain the following

\begin{prop}\label{gbocharprop} Let $\al$ and $\be$ with tracking chains $\alvec$ and $\bevec$, respectively, as above.
We have $\al\leo\be$ if and only if either $\alvec\subseteq\bevec$ or
there exists the $\klex$-increasing chain of index pairs $(i_1,j_1+1),\ldots,(i_s,j_s+1)\in\dom(\alvec)$ of maximal length $s\ge 1$ where $j_r\ge 1$ for $r=1,\ldots,s$,
such that $(i_1,j_1+1)\in\dom(\bevec)$, $\alvec_{\restriction_{i_1,j_1}}=\bevec_{\restriction_{i_1,j_1}}$, $\alcp{i_1,j_1+1}<\becp{i_1,j_1+1}$,
\[\alvec_{\restriction_{i_s,j_s+1}}\subseteq\alvec\subseteq\me(\alvec_{\restriction_{i_s,j_s+1}}),\]
and $\chi^\taucp{i_r,j_r}(\alcp{i_r,j_r+1})=1$ at least whenever $(i_r,j_r+1)\not=(n,m_n)$. 
Setting $\al_r:=\ov(\alvec_{\restriction_{i_r,j_r+1}})$ for $r=1,\ldots,s$ as well as $\al^+_r:=\ov(\me(\alvec_{\restriction_{i_r,j_r+1}}))$ 
for $r$ such that $\chi^\taucp{i_r,j_r}(\alcp{i_r,j_r+1})=1$ and $\al^+_s:=\al$ if $\chi^\taucp{i_s,j_s}(\alcp{i_s,j_s+1})=0$ we have
\[\al_1\ktwo\ldots\ktwo\al_s\letwo\al\leo\al^+_s\lo\ldots\lo\al^+_1\lo\ov(\bevec_{\restriction_{i_1,j_1+1}})\leo\be.\]
For $\be=\lh(\al)$ the cases $\alvec\subseteq\bevec$ and $s=1$ with $(i_1,j_1+1)=(n,m_n)$ correspond to the situation
$\gbo(\alvec)=(n,m_n)$, while otherwise we have $\gbo(\alvec)=(i_1,j_1+1)$.\qed
\end{prop}

\noindent{\bf Remark.}
Note that the above index pairs characterize the relevant sub-maximal $\nu$-indices in the initial chains of $\alvec$
with respect to $\bevec$ and leave the intermediate steps of maximal ($\me$-) extension along the iteration.
Using Lemma 5.5 of \cite{CWc} we observe that the sequence $\taucp{i_1,j_1},\ldots,\taucp{i_s,j_s}$ of bases in the above proposition satisfies
\begin{equation}\label{baseseqeq}
\taucp{i_1,j_1}<\ldots<\taucp{i_s,j_s}\quad\mbox{ and }\quad\taucp{i_s,j_s}<\taucp{i,1} \quad\mbox{ for every }i\in(i_s,n],
\end{equation}
so that in the case where $\al\lo\be$ and $\alvec\not\subseteq\bevec$ we have $\al\lo\gs(\al)\leo\be$ with $\taucp{i_s,j_s}\mid\rho_n\minusp1$.

\begin{lem}\label{subresplem}
The relation $\sub$ of initial chain on $\TC$ respects the ordering $\letc$ and hence also the characterization of $\leo$ on $\TC$.
\end{lem}
{\bf Proof.}
Suppose tracking chains $\alvec,\bevec,\gavec\in\TC$ satisfy $\alvec\letc\bevec\letc\gavec$ and $\alvec\sub\gavec$.
The case where any two chains are equal is trivial. We may therefore assume that $\ov(\alvec)<\ov(\bevec)<\ov(\gavec)$ as 
$(\TC,\letc)$ and $(\oneinf,<)$ are order-isomorphic. Since $\alvec\sub\gavec$ we have $\ov(\alvec)\lo\ov(\gavec)$, hence
$\ov(\alvec)\lo\ov(\bevec)$ as $\leo$ respects $\le$. Assume toward contradiction that $\alvec\not\sub\bevec$, 
whence by Proposition \ref{gbocharprop} there exists $(i,j+1)\in\dom(\alvec)\cap\dom(\bevec)$ such that 
$\alvec_{\restriction_{i,j}}=\bevec_{\restriction_{i,j}}$ and $\alcp{i,j+1}<\becp{i,j+1}$. 
But then for any $\devec\in\TC$ such that $\alvec\sub\devec$ we have 
\[\ov(\devec)<\ov(\alvecpl)\le\ov(\bevec)<\ov(\gavec),\] 
where the tracking chain $\alvecpl:=\alvec_{\restriction_{i,j+1}}[\alcp{i,j+1}+1]$
is the result of changing the terminal index $\alcp{i,j+1}$ of $\alvec_{\restriction_{i,j+1}}$ to $\alcp{i,j+1}+1$. 
\qed

We may now return to the issue of closure under $\lh$, which has a convenient sufficient condition on the basis of the following 
\begin{defi} A tracking chain $\alvec\in\TC$ is called \emph{convex} if and only if every $\nu$-index in $\alvec$
is maximal, i.e.\ given by the corresponding $\mu$-operator.
\end{defi}

\begin{cor}\label{lhclscor}
Let $\alvec\in\TC$ be convex and $M\finsub\TC$ be weakly spanning above $\alvec$.
Then $M$ is closed under $\lh$.
\end{cor}
{\bf Proof.}
This is a consequence of Proposition \ref{gbocharprop}, Corollary \ref{leoclosurecor}, Lemma \ref{meclosurelem}, 
and equation \ref{lhequation}.
\qed

While it is easy to observe that in $\Rtwo$ the relation $\leo$ is a forest that respects $\le$ and the relation $\letwo$ is a forest contained in $\leo$
which respects $\leo$, we can now conclude that this also holds for the arithmetical formulations of $\leo$ and $\letwo$ in $\Ctwo$, 
without referring to the results in Section 7 of \cite{CWc}.  

\begin{cor} Consider the arithmetical characterizations of $\leo$ and $\letwo$ on $\oneinf$.
The relation $\letwo$ respects $\leo$, i.e.\ whenever $\al\leo\be\leo\ga<\oneinf$ and $\al\letwo\ga$, then $\al\letwo\be$.
\end{cor}
{\bf Proof.} In the case $\bevec\subseteq\gavec$ this directly follows from Lemma \ref{letwosuclem}, 
while otherwise we additionally employ Proposition \ref{gbocharprop} and property \ref{baseseqeq}.
\qed

\section{Conclusion}
In the present article, the arithmetical characterization of the structure $\Ctwo$, which was established in Theorem 7.9 and 
Corollary 7.13 of \cite{CWc}, 
has been analyzed and shown to be elementary recursive. We have seen that finite isomorphism
types of $\Ctwo$ (in its arithmetical formulation) are contained in the class of ``respecting forests'', cf.\ \cite{C01}, over the language 
$(\le,\leo,\letwo)$. In a subsequent article \cite{W} we will establish the converse by providing an effective assignment of isominimal realizations 
in $\Ctwo$ to arbitrary respecting forests. We will provide an algorithm to find pattern notations for the ordinals in $\Ctwo$, and will
conclude that the union of isominimal realizations of respecting forests is indeed the core of $\Rtwo$, i.e.\ the structure $\Ctwo$ in its 
semantical formulation based on $\Sigma_i$-elementary substructures, $i=1,2$.
As a corollary we will see that the well-quasi orderedness of respecting forests with respect to coverings, which was shown by Carlson in \cite{C16},
implies (in a weak theory) transfinite induction up to the proof-theoretic ordinal $\oneinf$ of $\kplnod$.

We expect that the approaches taken here and in our treatment of the structure $\Ronepl$, see \cite{W07b} and \cite{W07c}, will naturally 
extend to an analysis of the structure $\Rtwopl$ and possibly to structures of patterns of higher order.
A subject of ongoing work is to verify that the core of $\Rtwopl$ matches the proof-theoretic strength of a limit of $\mathrm{KPI}$-models.  

\section*{Acknowledgements}
I would like to express my gratitude to Professor Ulf Skoglund for encouragement and support of my research.
I would like to thank Dr.\ Steven D.\ Aird for editing the manuscript.


\end{document}